\documentclass[english]{article}

\usepackage{babel}
\usepackage[latin1]{inputenc}
\usepackage{amsfonts}
\usepackage{amsmath}
\usepackage{mathrsfs}
\usepackage{graphicx}
\usepackage{caption}

\newtheorem{thm}{Theorem}
\newtheorem{cor}[thm]{Corollary}
\newtheorem{lem}[thm]{Lemma}
\newtheorem{dfn}[thm]{Definition}
\newtheorem{exa}[thm]{Example}
\newtheorem{rem}[thm]{Remark}
\numberwithin{thm}{section}

\begin{document}

\section*{Ville Turunen:\\
Time-frequency analysis
on groups}

\paragraph{Abstract:}
{\it
Phase-space analysis or time-frequency analysis
can be thought as Fourier analysis
simultaneously both in time and in frequency,
originating from signal processing and quantum mechanics.
On groups having unitary Fourier transform,
we introduce and study a natural family of time-frequency transforms,
and investigate the related pseudo-differential operators.
}

\section{Introduction}

Time-frequency analysis is a subfield of Fourier analysis.
It studies ``time'' dependent signals
(functions or distributions),
presenting them simultaneously both in ``time'' and in ``frequency'',
and consequently manipulating them as sharp as possible.
Traditionally, ``time'' and ``frequency''
refer to real variables,
where the Fourier integral transform is the essential tool.
In this text, we establish time-frequency analysis
on those locally compact groups that allow a unitary Fourier transform.

Time-frequency transforms in Cohen's class
present signals as joint time-frequency distributions,
which are linked to the pseudo-differential operators
for manipulating signals.
The time-frequency concepts apply to the phase-space analysis,
e.g. for position-momentum presentations
of wavefunctions in quantum mechanics.
One of Cohen's original motivating examples in \cite{Cohen1966}
was the deduction of the Born--Jordan phase-space transform,
stemming from the Born--Jordan quantization of Heisenberg's matrix mechanics
\cite{Heisenberg,BornJordan,BornHeisenbergJordan}.
Time-frequency analysis has been studied for $p$-adic numbers \cite{Haran},
on more general locally compact commutative groups \cite{Kutyniok},
and on certain classes of locally compact groups \cite{MantoiuRuzhansky}.
However, our treatise is not reduced to those works.

For a compact group, the time-frequency plane
is the Cartesian product of the group and its unitary dual.
Time-frequency transforms
will be ``time-frequency invariant''
sesquilinear mappings on pairs of test functions
(trigonometric polynomials, or Schwartz--Bruhat functions),
with values in the corresponding space of
matrix-valued test functions on the time-frequency plane.
In the non-commutative setting,
the ``frequency modulations'' require careful rethinking.
Euclidean time-frequency analysis is usually built around
the symmetric Wigner transform,
corresponding to the Weyl pseudo-differential quantization.
However, groups often lack suitable scalings,
so we build our time-frequency analysis
around the always existing Rihaczek or Kohn--Nirenberg transform:
this could have been the starting point for the Euclidean theory.
A time-frequency transform dictates a pseudo-differential quantization,
and we shall study this connection.
On compact Lie groups, the Kohn--Nirenberg quantization
has been treated e.g. in
\cite{Taylor,RuzhanskyTurunen,gardingRuzhanskyTurunen,RuzhanskyTurunenWirth,Fischer}.
The compact group results are finally generalized
to those locally compact groups that allow a unitary Fourier transform.

\section{On Euclidean time-frequency analysis}

To motivate our definitions for time-frequency analysis on compact groups $G$,
let us briefly explain how analogous concepts can be presented
on Euclidean spaces $\mathbb R^n$, avoiding technicalities.
The general background is presented in the monographs
\cite{Cohen1995} and \cite{Grochenig}.
To underline the similarities, we use quite similar notions
both on $G$ and on $\mathbb R^n$.
Signals are nice-enough functions $u:\mathbb R^n\to\mathbb C$.
We call variables $x,y\in\mathbb R^n$ {\it time-like} (or {\it position-like})
and variables $\xi,\eta\in\widehat{\mathbb R}^n\cong\mathbb R^n$
{\it frequency-like} (or {\it momentum-like}).
The starting point is formula
\begin{equation}
  u(x) = \iint {\rm e}^{{\rm i}2\pi(x-y)\cdot\eta}\,u(y)\,{\rm d}y\,{\rm d}\eta
\end{equation}
for the Schwartz test functions $u\in\mathscr S(\mathbb R^n)$.
Define the {\it Fourier transform} $\widehat{u}$ by
\begin{equation}
  \widehat{u}(\eta) :=
  \int {\rm e}^{-{\rm i}2\pi y\cdot\eta}\,u(y)\,{\rm d}y.
\end{equation}
From the Schwartz test function space $\mathscr S(\mathbb R^n)$,
the Fourier transform extends to a unitary operator
$\mathscr F:L^2(\mathbb R^n)\to L^2(\widehat{\mathbb R}^n)$:
in other words, $\mathscr F=(u\mapsto\widehat{u})$ is a linear bijection
satisfying
\begin{equation}
  \langle u,v\rangle := \langle \widehat{u},\widehat{v}\rangle,
\end{equation}
where Hilbert space $L^2(\mathbb R^n)$ has the inner product defined by
\begin{equation}
  \langle u,v\rangle = \int u(x)\,v(x)^\ast\,{\rm d}x,
\end{equation}
where $\lambda^\ast$ is
the complex conjugate of $\lambda\in\mathbb C$.
Signal $u$ has the {\it norm} $\|u\|=\langle u,u\rangle^{1/2}$ and
the {\it energy} $\|u\|^2=\langle u,u\rangle$.
The {\it symplectic Fourier transform} is then
$F=\mathscr F^{-1}\otimes\mathscr F$,
taking functions on the {\it time-frequency plane} (or {\it phase-space})
$\mathbb R^n\times\widehat{\mathbb R}^n$
to functions on the {\it ambiguity plane}
$\widehat{\mathbb R}^n\times\mathbb R^n$.
A Cohen class {\it time-frequency transform} $D$
of signals $u,v$
is $D(u,v):\mathbb R^n\times\widehat{\mathbb R}^n\to\mathbb C$,
\begin{eqnarray}
  D(u,v)(x,\eta)
    & = & F^{-1}\left(\phi\ FW(u,v) \right)(x,\eta) \\
    & = & \iint {\rm e}^{-{\rm i}2\pi y\cdot\eta}
          \,{\rm e}^{+{\rm i}2\pi x\cdot\xi}\,\phi(\xi,y)\,FW(u,v)(\xi,y)
          \,{\rm d}\xi\,{\rm d}y,\label{EQ:TFtransform}
\end{eqnarray}
where $\phi:\widehat{\mathbb R}^n\times\mathbb R^n\to\mathbb C$
is the {\it ambiguity kernel},
$W(u,v):\mathbb R^n\times\widehat{\mathbb R}^n\to\mathbb C$
is the {\it Wigner transform},
\begin{equation}\label{DEFN:Wigner}
  W(u,v)(x,\eta) := \int {\rm e}^{-{\rm i}2\pi y\cdot\eta}
  \,u(x+y/2)\,v(x-y/2)^\ast\,{\rm d}y,
\end{equation}
and $FW(u,v):\widehat{\mathbb R}^n\times\mathbb R^n\to\mathbb C$
is the {\it ambiguity transform},
\begin{equation}
  FW(u,v)(\xi,y) = \int {\rm e}^{-{\rm i}2\pi x\cdot\xi}
  \,u(x-y/2)\,v(x-y/2)^\ast\,{\rm d}x.
\end{equation}
As pointed out by Gr\"ochenig in \cite{Grochenig},
in the literature there is no precise definition of a Cohen class transform $D$.
Informally, such $D$ is obtained by smoothing the Wigner transform
by some tempered distribution as the convolution kernel.
In light of the time-frequency results in the sequel,
we would suggest that the ambiguity kernel $\phi$
should be a smooth function with polynomially bounded derivatives:
in other words, then we would have a Schwartz multiplier
$$
  \left(h\mapsto F^{-1}(\phi\,Fh)\right):
  \mathscr S(\mathbb R^n\times\widehat{\mathbb R}^n)
  \to\mathscr S(\mathbb R^n\times\widehat{\mathbb R}^n).
$$
Indeed, the literature examples of ambiguity kernels $\phi$ seem
to be smooth with polynomially bounded derivatives.
Moreover, those examples in the literature are typically bounded with
$|\phi(\xi,y)|\leq 1$,
which yields the $L^2$-boundedness
$$
  \|D(u,v)\|_{L^2} \leq \|u\|_{L^2} \|v\|_{L^2}.
$$
Hence $(u,v)\mapsto D(u,v)$ is sesquilinear:
$u\mapsto D(u,v)$ is linear, and $v\mapsto D(u,v)$ conjugate-linear.
The idea is that
the {\it time-frequency distribution}
$D[u]:=D(u,u)$ would be a quasi-energy density for signal $u$
(or a quasi-probability density for wavefunction $u$).
If $v(x)={\rm e}^{+{\rm i}2\pi x\cdot\xi}\,u(x-y)$ then
\begin{equation}
  D[v](x,\eta) = D[u](x-y,\eta-\xi),
\end{equation}
reflecting the idea that $v$ is
``$u$ shifted in time-frequency by $(y,\xi)$''.

For example, if the ambiguity kernel $\phi$ in \eqref{EQ:TFtransform}
is given by $\phi(\xi,y):={\rm e}^{{\rm i}2\pi(\xi\cdot y)\tau}$
for $\tau\in\mathbb R$, this defines the {\it Rihaczek-$\tau$-transform}
$D=R_\tau$, where
\begin{equation}
  R_\tau(u,v)(x,\eta) = \int {\rm e}^{-{\rm i}2\pi y\cdot\eta}
  \,u(x+(\tau+1/2)y)\,v(x+(\tau-1/2)y)^\ast\,{\rm d}y.
\end{equation}
Sometimes $W_\tau:=R_{\tau+1/2}$
is called the {\it Wigner-$\tau$}
or {\it Shubin-$\tau$ transform}.
Transforms $R_\tau$ and $R_{-\tau}$ are conjugates to each other
in the sense that
$R_\tau(u,v)(x,\eta)^\ast = R_{-\tau}(v,u)(x,\eta)$.
Especially, $R_0=W$, the Wigner transform.
The {\it Kohn--Nirenberg transform} (or the {\it Rihaczek transform})
is $R:=R_{-1/2}$,
which will be the starting point for time-frequency analysis on groups.
The {\it anti-Kohn--Nirenberg transform} refers to $R_{+1/2}$.
Here $D=R$ with
$\phi(\xi,y)={\rm e}^{-{\rm i}\pi\xi\cdot y}$, giving
\begin{equation}
 R(u,v)(x,\eta) = u(x)\,{\rm e}^{-{\rm i}2\pi x\cdot\eta}\,\widehat{v}(\eta)^\ast.
\end{equation}
It is easy to check that
$$
  \|R_\tau(u,v)\|_{L^2(\mathbb R^n\times\widehat{\mathbb R}^n)}
  = \|u\|_{L^2(\mathbb R^n)} \|v\|_{L^2(\mathbb R^n)}.
$$
Hence the {\it Born--Jordan transform} $Q$ defined by the integral average
\begin{equation}
  Q(u,v)
  = \int_0^1 W_\tau(u,v)\,{\rm d}\tau
  = \int_{-1/2}^{1/2} R_\tau(u,v)\,{\rm d}\tau
\end{equation}
satisfies
$$
  \|Q(u,v)\|_{L^2(\mathbb R^n\times\widehat{\mathbb R}^n)}
  \leq \|u\|_{L^2(\mathbb R^n)} \|v\|_{L^2(\mathbb R^n)},
$$
The ambiguity kernel of $D=Q$ satisfies
$\displaystyle\phi(\xi,y)
=\int_0^1 {\rm e}^{{\rm i}2\pi\xi\cdot y(\tau-1)}\,{\rm d}\tau
= {\rm sinc}(\xi\cdot y)$,
where ${\rm sinc}(t)=\sin(\pi t)/(\pi t)$ for $t\not=0$.

Let $\phi$ be the ambiguity kernel of time-frequency transform $D$.
Property $\phi(0,0)=1$ corresponds to the normalization
\begin{equation}\label{EQ:normalization}
  \iint D(u,v)(x,\eta)\,{\rm d}\eta\,{\rm d}x = \langle u,v\rangle.
\end{equation}
Properties $\phi(\xi,0)=1$ and $\phi(0,y)=1$
correspond respectively to the margins
\begin{equation}
  \int D(u,v)(x,\eta)\,{\rm d}\eta = u(x)\,v(x)^\ast,\quad
  \int D(u,v)(x,\eta)\,{\rm d}x = \widehat{u}(\eta)\,\widehat{v}(\eta)^\ast.
\end{equation}
Property $|\phi(\xi,y)|\equiv 1$ corresponds to
so-called {\it Moyal identity} \cite{Moyal}
\begin{equation}
  \langle D(u,v),D(f,g)\rangle = \langle u,f\rangle\,\langle v,g\rangle^\ast.
\end{equation}
In applied sciences and engineering,
perhaps the most common time-frequency transforms
date back to Gabor's work \cite{Gabor}:
these transforms $D$
are of the form
\begin{equation}
  D(u,v)(x,\eta) := \mathscr G_wu(x,\eta)\ \mathscr G_wv(x,\eta)^\ast,
\end{equation}
where the {\it $w$-windowed short-time Fourier transform} (STFT)
is defined by
\begin{equation}
  \mathscr G_w u(x,\eta) := \int {\rm e}^{-{\rm i}2\pi y\cdot\eta}
  \,u(y)\,w(y-x)^\ast\,{\rm d}y,
\end{equation}
where $\phi=FW[w]^\ast$.
Then the normalization \eqref{EQ:normalization}
means $\|w\|^2 = \langle w,w\rangle = 1$,
and then $D[u]=D(u,u)$ is called the {\it $w$-spectrogram} of $u$.

Once choosing a time-frequency transform $D$,
it defines the {\it $D$-quantization} $a\mapsto a^D$
by the $L^2$-duality
\begin{equation}
  \langle u,a^D v\rangle = \langle D(u,v),a\rangle.
\end{equation}
Here the weight function  $a:\mathbb R^n\times\widehat{\mathbb R}^n\to\mathbb C$
is called a {\it symbol} of {\it pseudo-differential operator}
$a^D=(v\mapsto a^D v)$.
Conversely, time-frequency transform $D$ can be recovered
from the quantization map $a\mapsto a^D$,
whose properties reflect the properties of $D$.
Wigner-$\tau$-transform $W_\tau=R_{\tau-1/2}$ corresponds to
so-called {\it Weyl-$\tau$-quantization} $a\mapsto a^{W_\tau}$,
\begin{equation}
  a^{W_\tau} v(x) = \iint {\rm e}^{{\rm i}2\pi(x-y)\cdot\eta}
  \,a(x+\tau(y-x),\eta)
  \, v(y)\,{\rm d}y\,{\rm d}\eta,
\end{equation}
Especially, the Wigner transform $W=W_{1/2}=R_0$
corresponds to the {\it Weyl quantization} $a\mapsto a^W$.
The Rihaczek (or Kohn--Nirenberg) transform $R=W_0=R_{-1/2}$ corresponds to the
{\it Kohn--Nirenberg quantization} $a\mapsto a^R$,
\begin{equation}
  a^R v(x) = \iint {\rm e}^{{\rm i}2\pi(x-y)\cdot\eta}\,a(x,\eta)
  \, v(y)\,{\rm d}y\,{\rm d}\eta
  = \int {\rm e}^{{\rm i}2\pi x\cdot\eta}\,a(x,\eta)\,\widehat{v}(\eta)
  \,{\rm d}\eta.
\end{equation}
The Born--Jordan quantization
$\displaystyle a\mapsto a^Q
= \int_0^1 a^{W_\tau}\,{\rm d}\tau
= \int_{-1/2}^{1/2} a^{R_\tau}\,{\rm d}\tau$
satisfies
\begin{equation}
  a^Q v(x) = \iint {\rm e}^{{\rm i}2\pi(x-y)\cdot\eta} \int_0^1 a(x+\tau(y-x),\eta)
  \,{\rm d}\tau\,v(y)\,{\rm d}y\,{\rm d}\eta.
\end{equation}
Weyl introduced his quantization in 1927 in \cite{Weyl},
and Wigner his distribution in 1932 in \cite{Wigner} for quantum mechanics.
The Wigner distribution was independently discovered in \cite{Ville},
with applications to signal processing.
The Born--Jordan quantization was implicit in \cite{BornJordan}
for polynomial symbols,
but in the modern sense the Born--Jordan distribution
was deduced by Cohen in \cite{Cohen1966}.
The Kohn--Nirenberg quantization arose from the studies
\cite{KohnNirenberg,Hormander}
by H\"ormander, Kohn and Nirenberg.

\section{Euclidean revision}

On a compact group $G$, we cannot expect to find a reasonable analogy to
Wigner transform $W$, which is the central object
in the Euclidean case presented above.
This is simply because analogies to the Euclidean scaling
$(y\mapsto y/2):\mathbb R^n\to\mathbb R^n$
are missing on a typical compact group $G$.
This problem does not disappear by a naive doubling change of variable
in the integral formula: see Example~\ref{EXA:falseWigner}.
Of course, on the odd-order cyclic group
$\mathbb Z/N\mathbb Z$,
such scalings $y\mapsto y/2$ exist in modular arithmetic.

On the other hand,
there is no necessity to start with the symmetric Wigner transform
in the Euclidean case, either.
Instead, we could have built the Cohen class theory
around the non-symmetric Kohn--Nirenberg quantization,
and this approach will work on compact groups, too.

There is another illuminating point of view:
Due to the time-frequency shift-invariance,
time-frequency transform $D$ is already encoded in data
\begin{equation}
  D(u,v)(0,0) = \langle D(u,v),\delta\rangle = \langle u,\delta^D v\rangle,
\end{equation}
where $\delta=\delta_{(0,0)}$ is the Dirac delta distribution
at the time-frequency origin $(0,0)\in\mathbb R^n\times\widehat{\mathbb R}^n$.
Despite such a highly singular symbol $\delta$,
pseudo-differential operator $\delta^D$
is typically rather well-behaving.
We call $\delta^D$ the {\it original localization operator},
as $\delta^D v$ tries to be the
``localization of $v$ to the time-frequency origin'',
which strictly speaking cannot be achieved
in view of the Heisenberg uncertainty principle.
If
\begin{equation}
  \delta^D v(z) = \int K_{\delta^ C}(z,y)\,v(y)\,{\rm d}y,
\end{equation}
i.e. if $K_{\delta^D}$ is the Schwartz distribution kernel of $\delta^D$,
then
\begin{equation}\label{EQ:kernelEuclidean}
  D(u,v)(x,\eta) = \iint u(x+z)\,{\rm e}^{-{\rm i}2\pi z\cdot\eta}
  \,K_{\delta^D}(z,y)^\ast\,{\rm e}^{+{\rm i}2\pi y\cdot\eta}\,v(x+y)^\ast
  \,{\rm d}z\,{\rm d}y.
\end{equation}
This formula suggests a natural variant for compact groups $G$,
where time shifts do not pose problems,
whereas frequency modulations are elusive.

\section{Fourier analysis on compact groups}

Let $e$ denote the neutral element of group $G$.
A topological group $G$ is a group and a Hausdorff space,
where the group operation $((x,y)\mapsto xy):G\times G\to G$
and the inversion $(x\mapsto x^{-1}):G\to G$
are continuous.

Time-frequency analysis on non-compact locally compact groups is
treated in Section~\ref{SEC:LCG},
generalizing most (but not all) of our compact case results.
In the sequel,
unless otherwise mentioned, $G$ is a {\it compact group}:
in other words, $G$ is a topological group with compact topology.
Monographs \cite{HewittRoss1} and \cite{HewittRoss2}
present background in Fourier analysis on compact groups.
From the Peter--Weyl theorem, which we shall review later,
it follows that such $G$
is isomorphic to a closed subgroup of the Cartesian product
of a family of unitary matrix groups.
If $G$ is commutative, instead of this multiplicative notation
for group operations,
it is common to use additive notation:
that is, instead of $xy,x^{-1},e$,
writing $x+y,-x,0$, respectively.

Let $C(G)$ be the vector space of continuous functions $u:G\to\mathbb C$,
endowed with the norm $u\mapsto\|u\|_{ C(G)}=\max\{|u(x)|:\, x\in G\}$.
Especially, the unit constant function ${\bf 1}=(x\mapsto 1):G\to\mathbb C$
belongs to $C(G)$.
Let
\begin{equation}
  \int u(x)\,{\rm d}x = \int_G u(x)\,{\rm d}x \in \mathbb C
\end{equation}
be the {\it Haar integral} of $u\in C(G)$:
the corresponding Haar measure is the unique
translation-invariant Borel probability measure on $G$.
We obtain the space $L^2(G)$ of square-integrable functions or {\it signals}
by completing $C(G)$ with respect to
the {\it norm} $\|u\|:=\langle u,u\rangle^{1/2}$
given by the the {\it inner product}
$(u,v)\mapsto\langle u,v\rangle$,
\begin{equation}
  \langle u,v\rangle := \int u(x)\,v(x)^\ast\,{\rm d}x.
\end{equation}
Here $\|u\|^2=\langle u,u\rangle$ is the {\it energy} of the signal.

A {\it unitary representation} of compact group $G$
on Hilbert space $\mathscr H_\eta$
is a strongly continuous group homomorphism
$\eta:G\to\mathscr U(\mathscr H_\eta)$
to the group $\mathscr U(\mathscr H_\eta)$
of unitary operators on $\mathscr H_\eta$.
Hence $\eta(xy)=\eta(x)\,\eta(y)$,
$\eta(x^{-1})=\eta(x)^{-1}=\eta(x)^\ast$,
$\eta(e)=I$ (the identity operator on $\mathscr H_{\eta}$).
The {\it Fourier coefficient} of $u\in L^2(G)$ at $\eta$
is the bounded linear operator
$\widehat{u}(\eta)=\mathscr F u(\eta):\mathscr H_\eta\to\mathscr H_\eta$
defined by
\begin{equation}
  \widehat{u}(\eta) = \mathscr F u(\eta) := \int u(x)\,\eta(x)^\ast\,{\rm d}x.
\end{equation}
The {\it left regular representation} of $G$
is $\pi_L:G\to\mathscr U(L^2(G))$ defined by
\begin{equation}
  \pi_L(y)u(x) := u(y^{-1}x)
\end{equation}
for almost all $x\in G$.
The left regular representation $\pi_L$ can be thought to embed the group $G$
into the ``rotations'' acting on Hilbert space $\mathscr H=L^2(G)$:
thus we can study the group
by tools of functional analysis.
Unitary representations $\xi,\eta$ of $G$ are {\it equivalent}
if there is a unitary isomorphism $U:\mathscr H_\xi\to\mathscr H_\eta$
such that
$$
  U\xi(x) = \eta(x) U
$$
for all $x\in G$.
The corresponding equivalence class is then denoted by $[\xi]=[\eta]$.
Unitary representation $\eta$ is called {\it irreducible}
if for operators $\eta(x)$ there are no non-trivial
simultaneous invariant subspaces of $\mathscr H_\eta$.
Let
$$
  \varepsilon=(x\mapsto 1):G\to\mathscr U(\mathbb C)
$$
denote the trivial irreducible unitary representation,
corresponding to ``zero frequency'', a unit signal with no oscillations.
We distinguist the trivial unitary representation $\varepsilon$
from the unit constant function
$
  {\bf 1}=(x\mapsto 1):G\to\mathbb C,
$
even though they are effectively the same.
This convention will clarify the treatise.

The {\it unitary dual} $\widehat{G}$ of $G$
consists of equivalence classes $[\eta]$
of irreducible unitary representations of $G$.
To make notation lighter,
instead of $[\eta]\in\widehat{G}$ we simply write
$\eta\in\widehat{G}$.
Due to the compactness of $G$,
for each $\eta\in\widehat{G}$,
Hilbert space $\mathscr H_\eta$ is finite-dimensional.
Hence in the sequel
we assume that $\eta(x)\in\mathbb C^{d_\eta\times d_\eta}$
is a unitary matrix of dimension $d_\eta\in\mathbb Z^+$:
there is such a choice in that equivalence class
$\eta\in\widehat{G}$.
The corresponding Fourier coefficient $\widehat{u}(\eta)$
is a matrix, belonging to $\mathbb C^{d_\eta\times d_\eta}$.
Function $u\in L^2(G)$ is called a {\it trigonometric polynomial}
if it has only finitely many non-zero Fourier coefficients:
in this sense,
trigonometric polynomials are band-limited signals.
Equivalently,
$u\in L^2(G)$ is a trigonometric polynomial
if and only if the span of $\{\pi_L(y)u:\ y\in G\}$
is a finite-dimensional vector space.
The space of trigonometric polynomials is denoted by $\mathscr T(G)$.

By the Peter--Weyl theorem, the left regular representation
can be decomposed to a direct sum of irreducible unitary representations
\begin{equation}
  \pi_L = \bigoplus_{\eta\in\widehat{G}} d_\eta\, \eta,
\end{equation}
corresponding to the Fourier decomposition of signals $u$:
in the sense of $L^2(G)$, there is the Fourier inverse formula
(Fourier series)
\begin{equation}
  u(x) = \sum_{\eta\in\widehat{G}} d_\eta\, {\rm tr}\left(
    \eta(x)\,\widehat{u}(\eta)\right),
\end{equation}
where ${\rm tr}$ is the usual matrix {\it trace}.
Here $\{\sqrt{d_\eta}\,\eta_{jk}:\ \eta\in\widehat{G},\ 1\leq j,k\leq d_\eta\}$
is an orthonormal basis for the Hilbert space $L^2(G)$.

Remember that ${\rm tr}(AB)={\rm tr}(BA)$,
but often ${\rm tr}(ABC)\not={\rm tr}(CBA)$.
In the sequel,
for matrix-valued functions $\widehat{a}$ on $\widehat{G}$,
we write ``non-commutative integrals''
\begin{equation}
  \int \widehat{a}(\eta)\,{\rm d}\eta
  :=
  \int_{\widehat{G}} {\rm tr}\left(\widehat{a}(\eta)\right)
    {\rm d}\mu_{\widehat{G}}(\eta)
  = \sum_{\eta\in\widehat{G}} d_\eta\,{\rm tr}(\widehat{a}(\eta))
\end{equation}
Here $\mu_{\widehat{G}}$ is the Plancherel measure.
We obtain
$$
  u(x) = \int \eta(x)\,\widehat{u}(\eta)\,{\rm d}\eta
  = \int \widehat{u}(\eta)\,\eta(x)\,{\rm d}\eta
  = \iint u(y)\,\eta(y^{-1}x)\,{\rm d}y\,{\rm d}\eta.
$$
Defining $\|\widehat{u}\|:=\langle\widehat{u},\widehat{u}\rangle^{1/2}$, where
\begin{equation}
  \langle\widehat{u},\widehat{v}\rangle
  := \int \widehat{u}(\eta)\,\widehat{v}(\eta)^\ast\,{\rm d}\eta,
\end{equation}
we obtain the Plancherel (or Parseval) identity
\begin{equation}\label{EQ:Plancherel}
   \langle u,v\rangle = \langle \widehat{u},\widehat{v}\rangle.
\end{equation}
Especially, $\|u\|^2=\|\widehat{u}\|^2$ is the conservation of energy.
Consequently, the Fourier transform $\mathscr F=(u\mapsto\widehat{u})$
is a Hilbert space isomorphism $\mathscr F:L^2(G)\to L^2(\widehat{G})$.
Furthermore, Fourier transform can also be vieved as linear isomorphisms
\begin{eqnarray}
  && \mathscr F=(u\mapsto\widehat{u}):\mathscr T(G)\to\mathscr T(\widehat{G}),\\
  && \mathscr F=(f\mapsto\widehat{f}):\mathscr T'(G)\to\mathscr T'(\widehat{G}),
\end{eqnarray}
where $\mathscr T'(G)$ is the space of {\it trigonometric distributions}
or formal trigonometric expansions $f$.
Here $\mathscr T'(\widehat{G})$ consists of all functions
$\widehat{f}$ on $\widehat{G}$
such that $\widehat{f}(\eta)\in\mathbb C^{d_\eta\times d_\eta}$
for each $\eta\in\widehat{G}$.
Elements
$\widehat{u}\in\mathscr T(\widehat{G})\subset\mathscr T'(\widehat{G})$
are those which have only finitely many non-zero Fourier coefficients.

On compact group $G$,
the algebra of test function can be enlarged from trigonometric $\mathscr T(G)$
to the Schwartz space (or Schwartz--Bruhat space) $\mathscr S(G)$,
introduced by Bruhat in \cite{Bruhat}.
Let $\mathscr J$ be the family of the closed normal subgroups $K$ of $G$
such that $G/K$ is isomorphic to a Lie group: for short, $G/K$ is a Lie group.
Endow $\mathscr J$ with the inverse inclusion order.
For $K\in\mathscr J$,
we identify
$u\in C^\infty(G/K)$ with $u\circ\pi_K:G\to\mathbb C$,
where $\pi_K=(x\mapsto xK):G\to G/K$ is the the quotient map.
Hence $C^\infty(G/K)\subset C(G)$.
The reflexive space of {\it Schwartz test functions}
is the inductive limit
$$
  {\mathscr S}(G) := \lim_{\longrightarrow} C^\infty(G/K)
$$
of the direct system
$\left( (C^\infty(G/K))_{K\in\mathscr J}, (f_{KL})_{K,L\in\mathscr J:\ K\subset L}
  \right)$,
where functions $f_{KL}:C^\infty(G/K)\to C^\infty(G/L)$
are defined by $f_{KL}(u)(xL):=u(xK)$.
The strong dual of the Schwartz space $\mathscr S(G)$
is the {\it Schwartz distribution space} $\mathscr S'(G)$,
and they are complete nuclear barreled spaces.

Function spaces are treated as subsets of distribution spaces,
and we have
$$
  \mathscr T(G)\subset\mathscr S(G)\subset C(G)
  \subset L^\infty(G)\subset L^2(G)\subset L^1(G)
  \subset\mathscr S'(G)\subset\mathscr T'(G).
$$
The Fourier transform can also be viewed as linear isomorphisms
\begin{eqnarray}
  && \mathscr F=(u\mapsto\widehat{u}):\mathscr S(G)\to\mathscr S(\widehat{G}),\\
  && \mathscr F=(f\mapsto\widehat{f}):\mathscr S'(G)\to\mathscr S'(\widehat{G}),
\end{eqnarray}
where $\mathscr S(\widehat{G})\subset L^2(\widehat{G})$
and $\mathscr S'(\widehat{G})\subset\mathscr T'(\widehat{G})$.

There is a positive central trigonometric approximate identity,
i.e. a net
of central positive trigonometric polynomials
$h_\alpha$ of unit $L^1$-norm such that
\begin{equation}
  \lim_\alpha \|u-h_\alpha\ast u\|_{L^1(G)} = 0
\end{equation}
for every $u\in L^1(G)$, see
\cite{HewittRoss2} (Theorem 28.53).
Let us present a brief related construction:
Let $\alpha=(U,m)$, where $m\in\mathbb Z^+$
and $U$ is a symmetric neighborhood of $e\in G$
meaning $U=U_eU_e$ for a neighborhood $U_e=U_e^{-1}$ of $e\in G$.
Choose central $f=f_U\in C(G)$ such that
$\|f\|_{L^2}=1$, and $f(x)=0$ whenever $x\not\in U$.
Approximate $f$ by central $g=g_{(U,m)}\in\mathscr T(G)$
such that $\|f-g\|_{C(G)} < 1/(m\|f\|_{C(G)})$.
Define central $h=h_{(U,m)}\in\mathscr T(G)$ by $h:=|g|^2/\|g\|_{L^2}^2$.
The index pairs $\alpha=(U,m)$ and $\beta=(V,n)$
have the partial order
\begin{eqnarray*}
  \alpha\leq\beta & \iff & V\subset U\ {\rm and}\ m\leq n.
\end{eqnarray*}
The functions $h_\alpha$
form a positive central trigonometric approximate identity.

Convolution $u\ast v$ of signals $u,v$ is the signal defined by
\begin{equation}
  u\ast v(x) := \int u(xy^{-1})\,v(y)\,{\rm d}y.
\end{equation}
Then $\widehat{u\ast v}=\widehat{v}\,\widehat{u}$, that is
$\widehat{u\ast v}(\eta)=\widehat{v}(\eta)\,\widehat{u}(\eta)$, as
$$
  \iint \eta(x)^\ast\, u(xy^{-1})\, v(y)\,{\rm d}y\,{\rm d}x
  = \int \eta(y)^\ast\,v(y) \int \eta(xy^{-1})^\ast\, u(xy^{-1})
  \,{\rm d}x\,{\rm d}y.
$$
The unitary dual $\widehat{G}$ does not have a group structure
when $G$ is non-commutative.
Nevertheless, we define a formal convolution by
\begin{equation}
  \widehat{u}\ast\widehat{v} :=
  \mathscr F\left((\mathscr F^{-1}\widehat{u})
    \,\mathscr F^{-1}\widehat{v}\right).
\end{equation}
Here we have commutativity
$\widehat{v}\ast\widehat{u}=\widehat{u}\ast\widehat{v}$
also on non-commutative groups $G$,
since multiplication of scalar-valued functions is commutative.

Matrix
$M=\begin{bmatrix} M_{jk}\end{bmatrix}\in\mathbb C^{d\times d}$
is {\it positive semi-definite} (or {\it positive}, for short) if
$$
  0\leq \langle Mz,z\rangle := \sum_{k=1}^d (Mz)_k\,z_k
  = \sum_{j,k=1}^d \overline{z_j}\,M_{jk}\,z_k.
$$
The Fourier series (or ``non-commutative integral'') over $\widehat{G}$
behaves much like the Haar integral over $G$.
For instance,
$$
  \int\widehat{u}(\eta)\,{\rm d}\eta=u(e),\quad\quad\quad
  \int u(x)\,{\rm d}x = \widehat{u}(\varepsilon).
$$
If $\widehat{u}\geq 0$
in the sense that $\widehat{u}(\eta)\geq 0$ for all $\eta\in\widehat G$
then $u(e)=\displaystyle\int\widehat{u}(\eta)\,{\rm d}\eta\geq 0$.

\begin{exa}
{\rm
For $1\leq p<\infty$
the {\it Schatten-$p$-norm} of a matrix $M\in\mathbb C^{d\times d}$ is
$$
  \|M\|_{S^p} := \left({\rm tr}(|M|^p)\right)^{1/p},
$$
where $|M|:=(M\,M^\ast)^{1/2}$.
The {\it operator norm} $\|M\|_{op}$ or the {\it Schatten-$\infty$-norm}
is the largest singular value of $M$,
$$
  \|M\|_{op} = \|M\|_{S^\infty} = \lim_{p\to\infty} \|M\|_{S^p},
$$
or alternatively
$\|M\|_{op} = \sup\{\|Mz\|_{\mathbb C^d}:\,\|z\|_{\mathbb C^d}\leq 1\}$,
where $\|z\|_{\mathbb C^d}^2 = \sum_{k=1}^d |z_k|^2$.
Here $\|M\|_{S_1}={\rm tr}(|M|)$ is the {\it trace class norm},
and $\|M\|_{HS}=\|M\|_{S^2}$
is the {\it Hilbert-Schmidt norm}.
The Lebesgue spaces $L^p(\widehat{G})$
have the norms given by
\begin{eqnarray}
  \|\widehat{u}\|_{L^p(\widehat{G})}
  & := & \left(\int |\widehat{u}(\eta)|^p\,{\rm d}\eta \right)^{1/p},\\
  \|\widehat{u}\|_{L^\infty(\widehat{G})}
  & := & \sup_{\eta\in\widehat{G}} \|\widehat{u}(\eta)\|_{op}.
\end{eqnarray}
If $1\leq p,q\leq\infty$ such that $1/p+1/q=1$, then
$$
  |u\ast v(e)|
  = \left| \int \widehat{v}(\eta)\,\widehat{u}(\eta)\,{\rm d}\eta \right|
  \leq \int |\widehat{v}(\eta)\,\widehat{u}(\eta)|\,{\rm d}\eta
  \leq \|\widehat{u}\|_{L^p(\widehat{G})}\,\|\widehat{v}\|_{L^q(\widehat{G})}.
$$
}
\end{exa}

\paragraph{Dirac and Kronecker deltas.}
In distributional sense,
the Fourier inverse formula
$\displaystyle
  u(x) = \int \eta(x)\ \widehat{u}(\eta)\,{\rm d}\eta
$
gives the expression
$$
  \delta_e(x) = \int \eta(x)\,{\rm d}\eta
$$
for the Dirac delta distribution $\delta_e\in C'(G)$ at $e\in G$.
Also, for $\eta\in\widehat{G}$,
$$
  \int \eta(x)^\ast\,{\rm d}x =
  \widehat{\bf 1}(\eta)
  = \delta_\varepsilon(\eta)\,I\in\mathbb C^{d_\eta\times d_\eta},
$$
where the Kronecker delta $\delta_\varepsilon$ at $\varepsilon\in\widehat{G}$
satisfies
$
  \delta_\varepsilon(\eta) :=
  \begin{cases} 1 & {\rm if}\ \varepsilon=\eta\in\widehat{G},\\
    0 & {\rm if}\ \varepsilon\not=\eta\in\widehat{G}.
    \end{cases}
$

\begin{exa}
{\rm
The compact commutative Lie groups $G$ are easy to list up to an isomorphism:
such a $G$ can be a product of a discrete cyclic group
$\mathbb Z/N\mathbb Z$
and a flat torus $\mathbb T^n=\mathbb R^n/\mathbb Z^n$
for some $N,n\in\mathbb N=\{0,1,2,3,\cdots\}$.
Let us review the notion above
in the familiar case of the torus group $G=\mathbb T^n=\mathbb R^n/\mathbb Z^n$.
The Haar measure on $G$ is given by the usual Lebesgue measure,
and for functions $u\in L^2(G)$
the traditional Fourier coefficient transform
$\widehat{u}:\mathbb Z^n\to\mathbb C$
is defined by
$$
  \widehat{u}(\eta) := \int_{\mathbb T^n} {\rm e}^{-{\rm i}2\pi y\cdot\eta}
  \,u(y)\,{\rm d}y.
$$
The inverse Fourier transform is given by the $L^2$-converging
Fourier series
$$
  u(x) = \sum_{\eta\in\mathbb Z^n} {\rm e}^{+{\rm i}2\pi x\cdot\eta}
  \,\widehat{u}(\eta).
$$
Here the irreducible unitary representations are one-dimensional
$$
  x\mapsto {\rm e}^{+{\rm i}2\pi x\cdot\eta},
$$
and we may obviously identify $\widehat{G}$
with $\mathbb Z^n$,
which is a non-compact discrete commutative group.
The convolutions are now given by
$$
  u\ast v(x) = \int_{\mathbb T^n} u(x-y)\,v(y)\,{\rm d}y,\quad\quad\quad
  \widehat{u}\ast\widehat{v}(\xi)
  = \sum_{\eta\in\mathbb Z^n} \widehat{u}(\xi-\eta)\,\widehat{v}(\eta).
$$
}
\end{exa}

\section{Hopf algebras of functions
and distributions}

Test function space $\mathscr T(G)$ of trigonometric polynomials
and $\mathscr S(G)$ of Schwartz functions
can be endowed with Hopf algebra structures.
Notice that
\begin{eqnarray*}
  \mathscr T(G\times G) & \cong & \mathscr T(G)\otimes\mathscr T(G),\\
  \mathscr S(G\times G) & \cong & \mathscr S(G) \hat{\otimes} \mathscr S(G),
\end{eqnarray*}
where $\otimes$ denotes the algebraic tensor product,
and $\hat{\otimes}$ the projective tensor product.
The commutative unital $C^\ast$-algebra $C(G)$ of continuous functions
has involution $\iota:C(G)\to C(G)$ given by $\iota u(x):=u(x)^\ast$.
Let us define mappings
\begin{eqnarray}
   m_0: C(G\times G)\to C(G),
  &&  m_0 w(x):=w(x,x),\\
  \eta_0:\mathbb C\to C(G),
  && \eta_0(\lambda):=\lambda\,{\bf 1},\\
  \Delta_0: C(G)\to C(G\times G),
  && \Delta_0 u(x,y):=u(xy),\\
  \varepsilon_0: C(G)\to\mathbb C,
  && \varepsilon_0 u(x) := u(e),\\
  S_0: C(G)\to C(G),
  && S_0 u(x) := u(x^{-1}).
\end{eqnarray}
When restricting these mappings
respectively to trigonometric polynomials to and Schwartz test functions,
$\mathscr T(G)$ and $\mathscr S(G)$ can be regarded as Hopf algebras.
By dualizing the structure of $\mathscr T(G)$, we obtain mappings
\begin{eqnarray*}
  && (m_1,\eta_1,\Delta_1,\varepsilon_1,S_1)\\
  & := & (\Delta_0',\varepsilon_0',m_0',\eta_0',S_0')
\end{eqnarray*}
where for $f,g\in\mathscr T'(G)$ we have
\begin{eqnarray}
  m_1(f\otimes g)=f\ast g, &&
     \mathscr F m_1(f\otimes g)(\xi)=\widehat{f}(\xi)\,\widehat{g}(\xi),\\
  \eta_1(\lambda) = \lambda\,\delta_e,
  && \mathscr F(\eta_1(\lambda))(\xi)=\lambda\,I
     \in\mathbb C^{d_\xi\times d_\xi},\\
  \Delta_1 f(x,y)=f(x)\,\delta_x(y),
  && \widehat{\Delta_1 f}(\xi\otimes\eta) = \widehat{f}(\xi\otimes\eta), \\
  \varepsilon_1(f) = \int f(x)\,{\rm d}x,
  && \varepsilon_1(f) = \widehat{f}(\varepsilon), \\
  S_1 f(x) = f(x^{-1}), && \widehat{S_1 f}(\eta)=\widehat{f}(\eta^\ast)^T.
\end{eqnarray}
Here $M^T$ is the transpose of matrix $M$, and
$\eta^\ast\in\widehat{G}$ is the contragredient representation
of $\eta\in\widehat{G}$, defined by
$\eta^\ast(x):=\eta(x^{-1})^T$.

\section{Symplectic Fourier transform}

We call $G\times\widehat{G}$ the {\it time-frequency plane}
(or the {\it position-momentum space}, or the {\it phase-space}),
where time-frequency points $(x,\eta)\in G\times\widehat{G}$
comprise of {\it time} $x\in G$
and of {\it frequency} $\eta\in\widehat{G}$.
We shall deal with Hilbert space $L^2(G\times\widehat{G})$,
where the inner product is given by
\begin{equation}
  \langle b,a\rangle
  = \iint b(x,\eta)\,a(x,\eta)^\ast
  \,{\rm d}\eta\,{\rm d}x.
\end{equation}
Here the matrix elements
of $x\mapsto a(x,\eta)\in\mathbb C^{d_\eta\times d_\eta}$
belong to $L^2(G)$
for all $\eta\in\widehat{G}$.
The {\it ambiguity plane}
\begin{equation}
  \widehat{G}\times G = \left\{(\xi,y):\ \xi\in\widehat{G},
    y\in G\right\}
\end{equation}
is the Fourier dual to the time-frequency plane $G\times\widehat{G}$
by the {\it symplectic Fourier transform} $F$, which is the linear isomorphism
\begin{equation}
  F=(\mathscr F\otimes I)(I\otimes\mathscr F^{-1}):
  L^2(G\times\widehat{G})\to L^2(\widehat{G}\times G).
\end{equation}
Thus if $a\in L^2(G\times\widehat{G})$
then $Fa\in L^2(\widehat{G}\times G)$,
\begin{equation}
  Fa(\xi,y) = \int \xi(x)^\ast \int \eta(y)\,a(x,\eta)
  \,{\rm d}\eta\,{\rm d}x.
\end{equation}
As in traditional signal processing,
here we may call $y\in G$ the {\it time-delay} or {\it lag} variable,
and $\xi\in\widehat{G}$ the {\it frequency-delay} or {\it Doppler} variable.
The inverse symplectic Fourier transform is then given by
\begin{equation}
  a(x,\eta) = \int \eta(y)^\ast \int \xi(x)\, Fa(\xi,y)
  \,{\rm d}\xi\,{\rm d}y.
\end{equation}
Then
$$
  \langle Fa,Fb \rangle = \langle a,b\rangle,\quad\quad\quad
  \|Fa\|^2 = \langle Fa,Fa\rangle = \langle a,a\rangle = \|a\|^2.
$$
Matrix-valued functions on $G\times\widehat{G}$ and $\widehat{G}\times G$
can be multiplied ``pointwise'':
$$
  (ab)(x,\eta):= a(x,\eta)\,b(x,\eta),\quad\quad\quad
  ((Fa)Fb)(\xi,y):=Fa(\xi,y)\,Fb(\xi,y).
$$
Then the {\it convolution} $a\ast b$ of $a,b$ on $G\times\widehat{G}$
is defined by
\begin{equation}
  a\ast b := F^{-1}((Fb)Fa).
\end{equation}
For example, $a\ast I = \lambda I$, where
\begin{equation}
  \lambda = Fa(\varepsilon,e)
  = \iint a(x,\eta)\,{\rm d}\eta\,{\rm d}x\in\mathbb C.
\end{equation}
We shall also need spaces of matrix-valued test functions
and distributions.
Especially,
we have linear isomorphisms
\begin{eqnarray}
  (I\otimes\mathscr F)
  & : & \mathscr S(G\times G)\to\mathscr S(G\times\widehat{G}), \\
  F & : & \mathscr S(G\times\widehat{G})\to\mathscr S(\widehat{G}\times G),
\end{eqnarray}
where we have the projective tensor product isomorphisms
\begin{eqnarray*}
  \mathscr S(G\times G)
  & \cong & \mathscr S(G)\hat{\otimes}\mathscr S(G), \\
   \mathscr S(G\times\widehat{G})
  & \cong & \mathscr S(G)\hat{\otimes}\mathscr S(\widehat{G}), \\
   \mathscr S(\widehat{G}\times G)
  & \cong & \mathscr S(\widehat{G})\hat{\otimes}\mathscr S(G).
\end{eqnarray*}
Then $\mathscr S'(\ldots)$ will denote the respective distribution space
corresponding to the test function space $\mathscr S(\ldots)$.

\section{Kohn--Nirenberg quantization}

The Kohn--Nirenberg quantization of pseudo-differential operators serves
as the starting point to acquire all the different time-frequency transforms.
The idea of the Kohn--Nirenberg pseudo-differential operators
on compact Lie groups was introduced by Taylor in \cite{Taylor},
and further investigated e.g. in
\cite{RuzhanskyTurunen,gardingRuzhanskyTurunen,RuzhanskyTurunenWirth,Fischer}.

\begin{dfn}
{\rm
The {\it Kohn--Nirenberg symbol} $a\in\mathscr S'(G\times\widehat{G})$
of linear mapping $B:\mathscr S(G)\to\mathscr S'(G)$
is defined by
\begin{equation}\label{DEFN:symbolKohnNirenberg}
  a(x,\eta) = \eta(x)^\ast B\eta(x),
\end{equation}
where matrix elements of $B\eta$ belong to $\mathscr S'(G)$. Then
\begin{equation}\label{EQ:operatorKohnNirenberg}
  Bv(x)
  = \int \eta(x)\,a(x,\eta)\,\widehat{v}(\eta)\,{\rm d}\eta
  = \int a(x,\eta)\,\widehat{v}(\eta)\,\eta(x)\,{\rm d}\eta,
\end{equation}
and we call $a^R:=B$ the {\it Kohn--Nirenberg pseudo-differential operator}
with symbol $a$.
The invertible mapping $a\mapsto a^R$ is called
the {\it Kohn--Nirenberg quantization}.
For $u,v\in\mathscr S(G)$, we define the corresponding
{\it Kohn--Nirenberg {\rm (or the {\it Rihaczek})} time-frequency transform}
$R(u,v)\in\mathscr S(G\times\widehat{G})$
by
\begin{equation}\label{EQ:quantizationKohnNirenberg}
  \langle u,a^{R}v\rangle = \langle R(u,v),a\rangle
\end{equation}
for all symbols $a\in\mathscr S'(G\times\widehat{G})$.
Then the {\it Kohn--Nirenberg ambiguity transform} is
$FR(u,v)=(\mathscr F\otimes I)(I\otimes\mathscr F^{-1})R(u,v)
\in\mathscr S(\widehat{G}\times G)$.
}
\end{dfn}

\begin{rem}
{\rm
Combining
\eqref{EQ:operatorKohnNirenberg} and
\eqref{EQ:quantizationKohnNirenberg},
we obtain
\begin{equation}\label{EQ:KN}
  R(u,v)(x,\eta) = u(x)\,\eta(x)^\ast\,\widehat{v}(\eta)^\ast
  \quad\in\quad\mathscr B(\mathscr H_\eta).
\end{equation}
Especially,
$R(u,v)(e,\varepsilon)=u(e)\,\widehat{v}(\varepsilon)^\ast\in\mathbb C$.
Notice that the same definition extends directly
to distributions $u,v\in\mathscr S'(G)$,
so that $R(u,v)\in\mathscr S'(G\times\widehat{G})$,
and then $FR(u,v)\in\mathscr S'(\widehat{G}\times G)$.
Moreover,
\begin{eqnarray}
  FR(u,v)(\xi,y)
  & = & \int \xi(x)^\ast \int \eta(y)\,R(u,v)(x,\eta)
  \,{\rm d}\eta\,{\rm d}x\\
  & = & \int \xi(x)^\ast\,u(x)\,v(xy^{-1})^\ast\,{\rm d}x
  \quad\in\quad \mathscr B(\mathscr H_\xi).
\end{eqnarray}
Especially,
$FR(u,v)(\varepsilon,e)
=\langle u,v\rangle
=\langle\widehat{u},\widehat{v}\rangle$.
Notice that $FR(u,v)(\xi,y) = \widehat{f_y\,}(\xi)$, where
$f_y(x)=u(x)\,v(xy^{-1})^\ast$. The Cauchy--Schwarz inequality yields
\begin{equation}
  \int |f_y(x)|\,{\rm d}x \leq \|u\|\,\|v\|.
\end{equation}
}
\end{rem}

\begin{rem}
{\rm
On a compact group $G$,
the Kohn--Nirenberg transform $R$ maps
$\mathscr T(G)\times\mathscr T(G)$ to $\mathscr T(G\times\widehat{G})$. Why?
Let $u,v\in\mathscr T(G)$. Since
$$
  u(x) = (I\otimes\varepsilon)\Delta u(x,y),\quad
  v(xy^{-1})^\ast = (I\otimes S)\Delta j v(x,y),
$$
this shows both $(x,y)\mapsto u(x)$ and $(x,y)\mapsto v(xy^{-1})^\ast$
belong to $\mathscr T(G\times G)$.
Hence also $(x,y)\mapsto u(x)\,v(xy^{-1})^\ast$
belongs to $\mathscr T(G\times G)$.
With a similar reasoning, we see that
the Kohn--Nirenberg transform maps
$\mathscr S(G)\times\mathscr S(G)$ to $\mathscr S(G\times\widehat{G})$.
}
\end{rem}

\section{Time-frequency transforms and quantizations}

Time-frequency transform $(u,v)\mapsto D(u,v)$
will be ``time-frequency invariant'',
taking test functions of time to
matrix-valued test functions of time-frequency. More precisely:

\begin{dfn}
{\rm
{\it Time-frequency transform}
$D:\mathscr S(G)\times\mathscr S(G)\to\mathscr S(G\times\widehat{G})$
is a mapping of the form
\begin{equation}\label{DEFN:TFtransform}
  D(u,v) := F^{-1}\left(\phi_D\,FR(u,v)\right),
\end{equation}
where $\phi_D\in\mathscr S'(\widehat{G}\times G)$
is the {\it ambiguity kernel} (or {\it Doppler-lag kernel}).
Time-frequency transform is called {\it band-limited}
if it maps $\mathscr T(G)\times\mathscr T(G)$
to $\mathscr T(G\times\widehat{G})$.
}
\end{dfn}

\begin{rem}
{\rm
Trigonometric function $u\in\mathscr T(G)$ is {\it band-limited}
in the sense that it has only finitely many non-zero Fourier coefficients.
For the Kohn--Nirenberg transform $R$,
notice that $\phi_R(\xi,y)=I$
for all $(\xi,y)\in\widehat{G}\times G$.
Thus we may have $\phi_D\not\in\mathscr S(\widehat{G}\times G)$.
Nevertheless,
\begin{equation}
  FD(u,{\bf 1})(\xi,y)=\phi_D(\xi,y)\,\widehat{u}(\xi)
\end{equation}
for all $u\in\mathscr T(G)$.
Hence the matrix elements of $y\mapsto\phi_D(\xi,y)$
belong to $\mathscr S(G)$ for each $\xi\in\widehat{G}$.
Band-limitedness of $D$ is equal
to that these matrix elements would be trigonometric polynomials:
for instance, the Kohn--Nirenberg transform $R$
is band-limited.
Time-frequency transform can also be expressed by
\begin{eqnarray}
  D(u,v)(x,\eta)
  & = & \int \eta(y)^\ast \int \xi(x)
  \,\phi_D(\xi,y)\,FR(u,v)(\xi,y)\,{\rm d}\xi\,{\rm d}y \\
  & = & R(u,v)\ast\psi_D(x,\eta),
\end{eqnarray}
where $\psi_D=F^{-1}(\phi_D)$
is the {\it time-frequency kernel} of $D$,
corresponding to the {\it ambiguity kernel} $\phi_D=F(\psi_D)$.
Sometimes we need the {\it time-lag kernel}
$\varphi_D = (I\otimes\mathscr F)\psi_D = (\mathscr F^{-1}\otimes I)\phi_D$.
Notice that the kernels
$$
  \psi_D(x,\eta),\quad \varphi_D(x,y),\quad \phi_D(\xi,y)
$$
contain the same information, with different variables
$x,y\in G$ and $\xi,\eta\in\widehat{G}$.
With the approach above, we have avoided finding ``frequency modulations''
on non-commutative groups;
the commutative case works still fine,
and yet we obtain many essential features also in the non-commutative setting.
}
\end{rem}

\begin{rem}
{\rm
If $u,v\in\mathscr T'(G)$
then $\phi_D\,FR(u,v)\in\mathscr T'(\widehat{G}\times G)$,
so that we can define
\begin{equation}
  D(u,v):=F^{-1}(\phi_D\,FR(u,v))\in\mathscr T'(G\times\widehat{G}).
\end{equation}
}
\end{rem}

\begin{dfn}
{\rm
Let $D$ be a time-frequency transform.
The corresponding {\it $D$-quantization} $a\mapsto a^D$ satisfies
\begin{equation}\label{DEFN:quantization}
  \langle u,a^D v\rangle = \langle D(u,v),a\rangle.
\end{equation}
Linear operators $a^D=(v\mapsto a^D v)$ are called
{\it $D$-pseudo-differential operators}.
}
\end{dfn}

In the sequel, we investigate how the properties of different kernels
affect the properties of the time-frequency transform $D$ and
the $D$-quantization $a\mapsto a^D$.
Due to \eqref{DEFN:quantization},
$a^Dv\in\mathscr S'(G)$ if $v\in\mathscr S(G)$ and
$a\in\mathscr S'(G\times\widehat{G})$.
Moreover,
if $v\in\mathscr S'(G)$ and $a\in\mathscr S(G\times\widehat{G})$,
then $a^Dv\in\mathscr S(G)$.
Thereby we have
\begin{eqnarray}
  a^D:\mathscr S(G)\to\mathscr S'(G)
  & {\rm if} & a\in\mathscr S'(G\times\widehat{G}),\\
  a^D:\mathscr S'(G)\to\mathscr S(G)
  & {\rm if} & a\in\mathscr S(G\times\widehat{G}).
\end{eqnarray}
Different quantizations can be linked
to the Kohn--Nirenberg case:

\begin{lem}\label{LEM:C2R}
Let $D$ be a time-frequency transform,
and let $a\in\mathscr S'(G\times\widehat{G})$.
Then $a^D=b^R$, where $Fb(\xi,y)=\phi_D(\xi,y)^\ast\,Fa(\xi,y)$.
\end{lem}

\paragraph{Proof.}
Noticing that
$$
  \langle D(u,v),a\rangle
  = \langle FD(u,v),Fa\rangle
  = \langle FR(u,v),Fb\rangle
  = \langle R(u,v),b\rangle,
$$
we obtain $\langle u,a^D v\rangle = \langle u,b^R v\rangle$.
\hfill {\bf QED}

\begin{dfn}
{\rm
For a time-frequency transform $(u,v)\mapsto D(u,v)$,
we call
\begin{equation}
  D[u]:=D(u,u)
\end{equation}
the {\it time-frequency distribution} of signal $u$.
Notice that $D[\lambda u]=|\lambda|^2 D[u]$ for all $\lambda\in\mathbb C$,
so define the equivalence class $[u]$ of {\it indistinguishable signals} by
\begin{equation}
  [u] := \left\{ \lambda u:\ \lambda\in\mathbb C,\ |\lambda|=1 \right\}.
\end{equation}
}
\end{dfn}

Value $D[u](x,\eta)\in\mathscr B(\mathscr H_\eta)$
presents an idealized operator-valued energy density
at time-frequency $(x,\eta)\in G\times\widehat{G}$
for a scalar-valued signal $u:G\to\mathbb C$.
With the complex scalars, numeric data families
\begin{equation}\label{EQ:uaCu}
  \left(\langle u,a^Du\rangle\right)_{u\in\mathscr S(G)}
  \quad{\rm and}\quad
  \left( \langle u,a^Dv\rangle\right)_{u,v\in\mathscr S(G)}
\end{equation}
mediate the same information.
Thereby the {\it invertibility} of time-frequency transform $D$
refers to the invertibility of the mapping $[u]\mapsto D[u]$.
This amounts to the properties of ambiguity kernel $\phi_D$.
Invertibility is not merely ``being bijective'',
it deals also with the numerical stability
(cf. the inverse problem for the traditional heat equation).
For invertibility, we need $\phi_D(\xi,y)$
to be invertible for almost every $(\xi,y)\in\widehat{G}\times G$,
and numerically that $\phi_D$ grows or decays at infinity at most polynomially.
The Kohn--Nirenberg transform is invertible, since
$$
  \int \eta(y)\,R[u](x,\eta)\,{\rm d}\eta
  = u(x)\,u(xy^{-1})^\ast.
$$

\begin{exa}
{\rm
An analogue of Wigner-$\tau$-pseudo-differential operators
on certain families of locally compact groups
was introduced and studied in \cite{MantoiuRuzhansky}.
On a compact group, this Wigner-$\tau$-quantization
would formally correspond to our time-frequency transform $D$,
which has the ambiguity kernel of the form
$$
  \phi_D(\xi,y) = \xi(\tau(y)),
$$
where $\tau:G\to G$ is a suitable function.
}
\end{exa}

\paragraph{Boundedness in energy.}
What if also $\phi_D\in L^\infty(\widehat{G}\times G)$?
In other words, $\phi_D$ would be bounded in the sense that
$\|\phi_D\|_{L^\infty} < \infty$ for
\begin{equation}
  \|\phi_D\|_{L^\infty}
  = \sup_{(\xi,y)\in\widehat{G}\times G} \|\phi_D(\xi,y)\|_{op},
\end{equation}
where $\|M\|_{op}$ is the spectral norm of operator $M$.
We obtain the following boundedness result on $L^2$-spaces,
where norms $\|f\|$ are the appropriate $L^2$-norms:

\begin{thm}\label{THM:L2boundedness}
Let
$\phi_D\in L^\infty(\widehat{G}\times G)$
for a time-frequency transform $D$.
Then
\begin{eqnarray}
  \|D(u,v)\| & \leq & \|\phi_D\|_{L^\infty}\,\|u\|\,\|v\|, \label{INEQ:L2C}\\
  \|a^Dv\| & \leq & \|\phi_D\|_{L^\infty}\,\|a\|\,\|v\|,
\end{eqnarray}
for all $u,v\in L^2(G)$ and $a\in L^2(G\times\widehat{G})$.
For the Kohn--Nirenberg transform,
$\|R(u,v)\|=\|u\|\,\|v\|$, and $\|a^Rv\|\leq \|a\|\,\|v\|$.
\end{thm}

\paragraph{Proof.}
In the special case of the Kohn--Nirenberg transform,
$\|\phi_R\|_{L^\infty}=1$
as $\phi_R(\xi,y)=I$ for all $(\xi,y)\in\widehat{G}\times G$.
Moreover,
\begin{eqnarray*}
  \|R(u,v)\|^2
  & = & \langle R(u,v),R(u,v)\rangle \\
  & = & \iint u(x)\,\eta(x)^\ast\,\widehat{v}(\eta)^\ast
        \,\widehat{v}(\eta)\,\eta(x)\,u(x)^\ast\,{\rm d}\eta\,{\rm d}x\\
  & = & \int |u(x)|^2\,{\rm d}x \int \widehat{v}(\eta)^\ast\,\widehat{v}(\eta)
        \,{\rm d}\eta\\
  & = & \|u\|^2\,\|v\|^2.
\end{eqnarray*}
The $L^2$-norm is preserved in the symplectic Fourier transform:
$$
  \|D(u,v)\| = \|FD(u,v)\| = \| \phi_D\,FR(u,v)\|.
$$
Let $\|M\|_{HS}=({\rm tr}(M M^\ast))^{1/2}$
denote the Hilbert--Schmidt norm.
Recall that
$\|MN\|_{HS}\leq \|M\|_{op}\|N\|_{HS}$. Thereby
\begin{eqnarray*}
  \|\phi_D\,FR(u,v)\|^2
  & = & \iint
        \|\phi_D(\xi,y)\,FR(u,v)(\xi,y)\|_{HS}^2
        \,{\rm d}\xi \,{\rm d}y \\
  & \leq & \iint \|\phi_D(\xi,y)\|_{op}^2
           \, \|FR(u,v)(\xi,y)\|_{HS}^2
           \,{\rm d}\xi \,{\rm d}y \\
  & \leq & \|\phi_D\|_{L^\infty}^2\,\|FR(u,v)\|^2
  \ = \ \|\phi_D\|_{L^\infty}^2\,\|R(u,v)\|^2.
\end{eqnarray*}
Inequality \eqref{INEQ:L2C} follows from this, because
$\|R(u,v)\|=\|u\|\,\|v\|$.
Hence by the Cauchy--Schwarz inequality we obtain
\begin{eqnarray*}
  \left|\langle u,a^D v\rangle\right| = \left|\langle D(u,v),a\rangle\right|
  \leq \|D(u,v)\|\,\|a\|
  \leq \|\phi_D\|_{L^\infty}\,\|u\|\,\|v\|\,\|a\|,
\end{eqnarray*}
completing the proof.
\hfill {\bf QED}

${}$

Let us emphasize the invariance under the time translations,
in the sense that
$D[v](x,\eta)=D[u](yx,\eta)$ if $v(x)=u(yx)$.
The frequency modulations are more elusive in the non-commutative case,
but nevertheless, the message of the next result
is that the information can be ``shifted'' to
the specific point $(e,\varepsilon)$
in the time-frequency plane:

\begin{thm}\label{THM:localizationoriginal}
Time-frequency transform $D$ can be recovered
from the evaluation mapping
$(u\mapsto D[u](e,\varepsilon)):\mathscr S(G)\to\mathbb C$.
\end{thm}

\paragraph{Proof.}
For $u,v\in\mathscr S(G)$ we have
\begin{eqnarray*}
  D(u,v)(x,\eta)
  & = & \int \eta(y)^\ast \int \xi(x)\,\phi_D(\xi,y)
        \int \xi(t)^\ast\,u(t)\,v(ty^{-1})^\ast
        \,{\rm d}t\,{\rm d}\xi\,{\rm d}y \\
  & = & \int \eta(y)^\ast \int \varphi_D(t^{-1}x,y)\,u(t)\,v(ty^{-1})^\ast
        \,{\rm d}t\,{\rm d}y.
\end{eqnarray*}
Especially,
\begin{equation}\label{EQ:original}
  \langle u,\delta^D u\rangle = D[u](e,\varepsilon)
  = \int u(x) \left(\int \varphi_D(x^{-1},y^{-1}x)^\ast\,u(y)\,{\rm d}y
        \right)^\ast {\rm d}x,
\end{equation}
where $\delta=\delta_{(e,\varepsilon)}$
is the Dirac--Kronecker delta distribution at
$(e,\varepsilon)\in G\times\widehat{G}$.
Hence from knowing all $D[u](e,\varepsilon)$ we obtain $\varphi_D$
and thereby $D$.
\hfill {\bf QED}

\begin{rem}{\rm
Notice that in the statement of the previous Theorem on a compact group $G$,
we could replace the test function space $\mathscr S(G)$ by $\mathscr T(G)$.}
\end{rem}

\paragraph{Original localizations.}
Let us call $(e,\varepsilon)\in G\times\widehat{G}$
the {\it origin of the time-frequency plane} $G\times\widehat{G}$.
As seen in the proof of Theorem~\ref{THM:localizationoriginal} above,
the {\it original localization} pseudo-differential operator
$\delta^D=(\delta_{(e,\varepsilon)})^D:\mathscr S(G)\to\mathscr S'(G)$
encodes all the information about the time-frequency transform $D$.
The original localization $\delta^D$ is bounded on $L^2(G)$
if and only if
\begin{equation}\label{DEFN:uniformbound}
  |D(u,v)(e,\varepsilon)| = |\langle u,\delta^D v\rangle| \leq c\,\|u\|\,\|v\|
\end{equation}
for all $u,v\in\mathscr S(G)$, where $c<\infty$ is a constant.
Original localizations provide an alternative way to understand
time-frequency transforms.
Notice that if $K_{\delta^D}$ is the Schwartz integral kernel of
the original localization $\delta^D$, then
\begin{equation}\label{EQ:Cfrom0}
  D(u,v)(x,\eta) = \iint u(xz)\,\eta(z)^\ast
  \,K_{\delta^D}(z,y)^\ast\,\eta(y)\,v(xy)^\ast
  \,{\rm d}z\,{\rm d}y,
\end{equation}
in analogy to the Euclidean case \eqref{EQ:kernelEuclidean}.
Here $K_{\delta^D}(z,y)^\ast = \varphi_D(z^{-1},y^{-1}z)$,
that is $\varphi_D(x,y)=K_{\delta^D}(x^{-1},x^{-1}y^{-1})^\ast$.
Hence
\begin{equation}\label{EQ:ambiguitySchwartz}
  \phi_D(\xi,y) = \int K_{\delta^D}(x,xy^{-1})^\ast\,\xi(x)\,{\rm d}x.
\end{equation}
Moreover, if we define {\it amplitude} $a_D$ by
\begin{equation}
  a_D(x,y,\eta):=\int a(t,\eta)\,K_{\delta^D}(t^{-1}x,t^{-1}y)\,{\rm d}t
\end{equation}
then $\displaystyle a^D v(x) = \int K_{a^D}(x,y)\,v(y)\,{\rm d}y$
for the Schwartz kernel $K_{a^D}$:
\begin{equation}
  K_{a^D}(x,y) = \int \eta(y^{-1}x)\,a_D(x,y,\eta)\,{\rm d}\eta.
\end{equation}

\begin{exa}
{\rm
Since
$R(u,v)(e,\varepsilon)=u(e)\,\widehat{v}(\varepsilon)^\ast$,
the Kohn--Nirenberg original localization is given by
\begin{equation}
  \delta^{R}v(x)
  = \widehat{v}(\varepsilon)\,\delta_e(x)
  =\int v(y)\,{\rm d}y\ \delta_e(x).
\end{equation}
Here
$\delta^{R}:\mathscr S(G)\to\mathscr S'(G)$ is unbounded on $L^2(G)$
unless $G$ is finite.
Amplitude $a_R$ of $a^R$ satisfies $a_R(x,y,\eta)=a(x,\eta)$.
The so-called {\it anti-Kohn--Nirenberg transform} $R^\ast$ satisfies
$R^\ast(u,v)(e,\varepsilon)=\widehat{u}(\varepsilon)\,v(e)^\ast$.
Its original localization satisfies
$
  \delta^{(R^\ast)} v(x) = v(e),
$
and its amplitudes are given by
$a_{R^\ast}(x,y,\eta)=a(y,\eta)$.
}
\end{exa}

\section{Uncertainty and original localizations}

In this section we discuss original localizations
related to the Heisenberg uncertainty principle in quantum mechanics.
Our quantum states $u$ are unit vectors
in the Hilbert space $\mathscr H=L^2(G)$,
identifying states $u,v$ whenever $[u]=[v]$.
Bounded observables are self-adjoint operators
$A:\mathscr H\to \mathscr H$, and
\begin{equation}
  \begin{cases}
  \mu = \mu_A^u := \langle Au,u\rangle, & \\
  \sigma = \sigma_A^u := \|Au-\mu u\| &
  \end{cases}
\end{equation}
are the {\it expectation} and the {\it deviation}
of measurement $A$ in state $u$, respectively.
For instance, let $A=\sum_{\alpha\in J} \alpha\,P_\alpha$
where $P_\alpha$ is an orthogonal projection,
with distinct measured values $\alpha\in J\subset\mathbb R$.
Then the interpretation is the following:
in initial state $u$,
our measurement gives value $\alpha\in J$
with probability $\|P_\alpha u\|^2$,
and then $u$ collapses to state
$P_\alpha u/\|P_\alpha u\|$.
Let $A,B$ be bounded observables.
The {\it uncertainty observable} of the pair $(A,B)$ is
\begin{equation}
  -{\rm i}\hbar^{-1} [A,B] = -{\rm i}\hbar^{-1} (AB-BA),
\end{equation}
where we normalize the Dirac--Planck constant so that
$\hbar := (2\pi)^{-1}$. Applying Cauchy--Schwarz inequality, we obtain
the Heisenberg uncertainty inequality
\begin{equation}
  \left|\mu_{-{\rm i}\hbar^{-1}[A,B]}^u\right|
  \leq 2\hbar^{-1}\,\sigma_A^u\,\sigma_B^u.
\end{equation}

Suppose above $A$ would be a ``position operator''
and $B$ a ``momentum operator'':
$Au=fu$ and $\widehat{Bu}=\widehat{u}\,\widehat{g}$ (that is $Bu=g\ast u$),
initially with $f,g\in\mathscr S(G)$
(later considering $f,g\in\mathscr S'(G)$),
where for self-adjointness we should have
real-valued ``coordinate function'' $f$,
and $g(z)^\ast=g(z^{-1})$.
If $A,B$ are able to distinguish $(e,\varepsilon)\in G\times\widehat{G}$
in a reasonable fashion,
then a good candidate for an original localization $\delta^D$
would be given by
\begin{equation}
  \delta^D := -{\rm i}2\pi [A,B].
\end{equation}
Then
\begin{eqnarray*}
  D(u,v)(e,\varepsilon)
  & = & \iint {\rm i}2\pi\left(f(x)-f(y)\right) g(xy^{-1})^\ast\,u(x)\,v(y)^\ast
        \,{\rm d}y\,{\rm d}x.
\end{eqnarray*}
We shall return to this uncertainty commutator approach
when dealing with cyclic groups in Section~\ref{SEC:cyclic}.
If here $f=\delta_e\in\mathscr S'(G)$ and $g={\bf 1}$
then
$$
  D(u,v) = {\rm i}2\pi \left( R(u,v)-R^\ast(u,v) \right),
$$
where the conjugate transforms $R^\ast(u,v):=R(v,u)^\ast$
will be studied in Section~\ref{SEC:symmetry}.

\section{Symmetry}\label{SEC:symmetry}

\begin{dfn}
{\rm
Let $D:\mathscr S(G)\times\mathscr S(G)\to\mathscr S(G\times\widehat{G})$
be a time-frequency transform.
We define its {\it conjugate}
\begin{equation}\label{DEFN:conjugateC}
  D^\ast:\mathscr S(G)\times\mathscr S(G)\to\mathscr S(G\times\widehat{G})
\end{equation}
by $D^\ast(u,v):=D(v,u)^\ast$,
more precisely
\begin{equation}\label{DEFN:cC}
  D^\ast(u,v)(x,\eta) := D(v,u)(x,\eta)^\ast
\end{equation}
for all $u,v\in\mathscr S(G)$ and $(x,\eta)\in G\times\widehat{G}$.
We call time-frequency transform $D$ {\it symmetric}
if $D^\ast=D$.
The $D$-quantization is {\it symmetric}
if for all $u,v\in\mathscr S(G)$
$$
  \langle a^D u,v\rangle = \langle u,a^D v\rangle
$$
whenever $a\in\mathscr S(G\times\widehat{G})$
satisfies $a(x,\eta)^\ast=a(x,\eta)$ for all $(x,\eta)\in G\times\widehat{G}$.
}
\end{dfn}

\begin{thm}
Mapping $D^\ast$ defined in \eqref{DEFN:conjugateC},\eqref{DEFN:cC}
is a time-frequency transform.
Moreover, the following conditions are equivalent:
\begin{itemize}
\item[{\rm (a)}]
For all $u\in\mathscr S(G)$ we have $D[u](e,\varepsilon)\in\mathbb R$.
\item[{\rm (b)}]
Time-frequency transform $D$ is symmetric.
\item[{\rm (c)}]
The $D$-quantization is symmetric.
\item[{\rm (d)}]
The time-lag kernel
$\varphi_D=(I\otimes\mathscr F^{-1})\psi_D=(\mathscr F^{-1}\otimes I)\phi_D$
satisfies
\begin{equation}\label{EQ:symmetry}
  \varphi_D(x,y)^\ast = \varphi_D(yx,y^{-1}).
\end{equation}
\end{itemize}
\end{thm}

\begin{rem}
{\rm
The superficial non-symmetry in the appearance
of \eqref{EQ:symmetry}
is just due to the fact that the Kohn--Nirenberg transform
itself is not symmetric.
Moreover, in the statement of the previous Theorem on a compact group $G$,
we can replace the test function space $\mathscr S(G)$ by $\mathscr T(G)$.
}
\end{rem}

\paragraph{Proof.}
On the one hand,
\begin{eqnarray*}
  D(u,v)(x,\eta)
  & = &  \int \eta(y)^\ast \int \xi(x)\,\phi_D(\xi,y)\,FR(u,v)(\xi,y)
        \,{\rm d}\xi\,{\rm d}y \\
  & = &  \int \eta(y)^\ast \int \xi(x)\,\phi_D(\xi,y)
        \int \xi(z)^\ast\,u(z)\,v(zy^{-1})^\ast
        \,{\rm d}z\,{\rm d}\xi\,{\rm d}y \\
  & = &  \int \eta(y)^\ast \int \varphi_D(z^{-1}x,y)\,u(z)\,v(zy^{-1})^\ast
        \,{\rm d}z\,{\rm d}y.
\end{eqnarray*}
On the other hand,
\begin{eqnarray*}
  D(v,u)(x,\eta)^\ast
  & = & \int \eta(y) \int \varphi_D(z^{-1}x,y)^\ast\,v(z)^\ast\,u(zy^{-1})
        \,{\rm d}z\,{\rm d}y \\
  & = & \int \eta(y)^\ast \int \varphi_D(z^{-1}x,y^{-1})^\ast\,u(zy)
        \,v(z)^\ast\,{\rm d}z\,{\rm d}y \\
  & = & \int \eta(y)^\ast \int \varphi_D(yz^{-1}x,y^{-1})^\ast\,u(z)
        \,v(zy^{-1})^\ast\,{\rm d}z\,{\rm d}y,
\end{eqnarray*}
showing that
\begin{equation}\label{EQ:timelagconjugate}
  \varphi_{D^\ast}(x,y) = \varphi_D(yx,y^{-1})^\ast,
\end{equation}
and leading to the equivalence of conditions (b) and (d).
In the special case of $(x,\eta)=(e,\varepsilon)$ and $v=u$,
this gives also the equivalence of (d) and (a).
Moreover, if $D$ is symmetric and $a^\ast=a$, then
$$
  \langle a^D u,u\rangle
  = \langle u,a^D u\rangle^\ast
  = \langle D[u],a\rangle^\ast
  = \langle D[u]^\ast,a^\ast\rangle
  = \langle D[u],a\rangle
  = \langle u,a^D u\rangle,
$$
so that $a\mapsto a^D$ is also symmetric.
Thus (b) implies (c).
Now suppose $a\mapsto a^D$ is symmetric.
Let $(h_\alpha)_\alpha$ be a bounded left approximate identity
with $0\leq h_\alpha\in\mathscr S(G)$.
Define $a_\alpha\in\mathscr S(G\times\widehat{G})$ by
$a_\alpha(x,\eta):=h_\alpha(x)\,\delta_\varepsilon(\eta)I$.
Then $a_\alpha(x,\eta)^\ast=a_\alpha(x,\eta)$,
and
$$
  D[u](e,\varepsilon)
  = \langle D[u],\delta_{(e,\varepsilon)}\rangle
  = \lim_{\alpha}\langle D(u,u),a_\alpha\rangle
  = \lim_{\alpha}\langle u,(a_\alpha)^D u\rangle,
$$
which is real-valued due to the symmetry of the quantization.
Hence condition (c) implies (a).
\hfill {\bf QED}

\begin{rem}
{\rm
Clearly, $(D^\ast)^\ast = D$. Notice also that
\begin{equation}
  (a^D)^\ast = (a^\ast)^{(D^\ast)},
\end{equation}
and especially $(\delta^D)^\ast=\delta^{(D^\ast)}$.
This follows from
$$
  \langle u,(a^D)^\ast v\rangle
  = \langle v,a^D u\rangle^\ast
  = \langle D(v,u)^\ast,a^\ast\rangle
  = \langle D^\ast(u,v),a^\ast\rangle
  = \langle u,(a^\ast)^{(D^\ast)}v\rangle.
$$
}
\end{rem}

\begin{exa}
{\rm
The conjugate $R^\ast$ of the Kohn--Nirenberg transform $R$
satisfies
\begin{equation}
  R^\ast(u,v)(x,\eta) = R(v,u)(x,\eta)^\ast
  = \widehat{u}(\eta)\,\eta(x)\,v(x)^\ast.
\end{equation}
The corresponding pseudo-differential quantization satisfies
\begin{eqnarray*}
  \langle u,a^{(R^\ast)} v\rangle
  & = & \langle R^\ast(u,v),a\rangle \\
  & = & \iiint \widehat{u}(\eta)\,\eta(x)\,v(x)^\ast\,a(x,\eta)^\ast
        \,{\rm d}\eta\,{\rm d}x \\
  & = & \int u(y) \left( \iint a(x,\eta)\,v(x)\,\eta(x^{-1}y)
        \,{\rm d}\eta\,{\rm d}x\right)^\ast {\rm d}y,
\end{eqnarray*}
leading to
\begin{equation}
  a^{(R^\ast)} v(x)
  = \iint \eta(y^{-1}x)\,a(y,\eta)\,v(y)\,{\rm d}\eta\,{\rm d}y.
\end{equation}
Mapping $a\mapsto a^{(R^\ast)}$
is called the {\it anti-Kohn--Nirenberg quantization}.
It is easy to find that $\phi_{R^\ast}(\xi,y)=\xi(y)$.
}
\end{exa}

\begin{exa}
{\rm
Let $D$ be a time-frequency transform.
Then
$$
  D=\frac{D+D^\ast}{2}+{\rm i}\,\frac{D-D^\ast}{2{\rm i}},
$$
where the symmetric time-frequency transforms
$(D+D^\ast)/2$ and
$-{\rm i}(D-D^\ast)/2$ could be called the respective
{\it real} and {\it imaginary parts} of $D$.
}
\end{exa}

\begin{exa}\label{EXA:falseWigner}
{\rm
For the moment, let us try to introduce Wigner distribution
on compact groups $G$.
The Euclidean space Wigner transform \eqref{DEFN:Wigner} satisfies
\begin{eqnarray*}
  W(u,v)(x,\eta)
  & = & \int_{\mathbb R^n} {\rm e}^{-{\rm i}2\pi y\cdot\eta}
        \,u(x+y/2)\,v(x-y/2)^\ast\,{\rm d}y \\
  & = & 2^n \int_{\mathbb R^n} {\rm e}^{-{\rm i}2\pi 2z\cdot\eta}
        \,u(x+z)\,v(x-z)^\ast\,{\rm d}z.
\end{eqnarray*}
It would be tempting to define the ``Wigner transform ${\bf W}$''
of $u,v\in\mathscr S(G)$ by
\begin{equation}\label{DEFN:falseWigner}
  {\bf W}(u,v)(x,\eta) := \int \eta(z)^\ast\,u(xz)\,v(xz^{-1})^\ast\,{\rm d}z
\end{equation}
possibly up to a constant multiple, depending on $G$.
The problem here is that $(xz^{-1})^{-1}(xz)=z^2$
is not the lag $z$ in time.
Transform ${\bf W}$ would also be formally symmetric, as
${\bf W}(v,u)(x,\eta)^\ast={\bf W}(u,v)(x,\eta)$.
Nevertheless,
\begin{eqnarray*}
  F{\bf W}(u,{\bf 1})(\xi,y)
  & = & \int \xi(x)^\ast \int \eta(y) \int \eta(z)^\ast\,u(xz)
        \,{\rm d}z\,{\rm d}\eta\,{\rm d}x \\
  & = & \int \xi(x)^\ast u(xy)\,{\rm d}x \\
  & = & \xi(y)\,\widehat{u}(\xi).
\end{eqnarray*}
So for such ${\bf W}$ to be a time-frequency transform,
we would have $\phi_{\bf W}(\xi,y)=\xi(y)$,
meaning that ${\bf W}=R^\ast$, the anti-Kohn--Nirenberg transform:
this is possible only when $G=\{e\}$ is the trivial group of one element.
Hence, $\bf W$ defined in formula \eqref{DEFN:falseWigner}
is a dead-end in time-frequency analysis,
and it does not make sense to talk about a corresponding
Weyl-like pseudo-differential quantization:
especially, $(u,v)\mapsto {\bf W}(u,v)$
would not be modulation-invariant
for commutative $G\not=\{e\}$.
However, consider such a compact group $G$,
where $(y\mapsto y^2):G\to G$ is a bijection:
its inverse $(y\mapsto y^{1/2}):G\to G$ is a homeomorphism
of the compact Hausdorff space $G$.
Then
\begin{equation}
  W(u,v)(x,\eta) := \int \eta(y)^\ast\,u(xy^{1/2})\,v(xy^{-1/2})^\ast\,{\rm d}y
\end{equation}
defines the natural {\it Wigner transform} on $G$,
where $y^{-1/2}=(y^{1/2})^{-1}$,
and $\phi_W(\xi,y)=\xi(y^{1/2})$.
Especially, it is possible to define the Wigner time-frequency transform
on finite cyclic groups of odd order,
or on $p$-adic groups for primes $p\not=2$.
Related questions on commutative locally compact groups
have been treated in \cite{Kutyniok}.
}
\end{exa}

\section{Normalization, and time-frequency margins}

\begin{dfn}
{\rm
We call time-frequency transform $D$ {\it normalized} if
\begin{equation}\label{EQ:marginenergy}
  \iint D(u,v)(x,\eta)\,{\rm d}\eta\,{\rm d}x
  = \langle u,v\rangle = \langle \widehat{u},\widehat{v}\rangle
\end{equation}
for all $u,v\in\mathscr S(G)$.
Especially for $v=u$ formula \eqref{EQ:marginenergy}
yields the energy $\|u\|^2=\|\widehat{u}\|^2$.
We say that the $D$-quantization has {\it correct traces} if
\begin{equation}\label{EQ:trace}
  {\rm tr}(a^D) = \iint a(t,\eta)\,{\rm d}\eta\,{\rm d}t
\end{equation}
for all $a\in\mathscr S(G\times\widehat{G})$.
}
\end{dfn}

\begin{thm}
The following conditions are equivalent:
\begin{itemize}
\item[{\rm (a)}] $\iint D[{\bf 1}](x,\eta)\,{\rm d}\eta\,{\rm d}x = 1$.
\item[{\rm (b)}] Time-frequency transform $D$ is normalized.
\item[{\rm (c)}] The $D$-quantization has correct traces.
\item[{\rm (d)}]
The ambiguity kernel satisfies
$\phi_D(\varepsilon,e)=1\in\mathbb C$.
\end{itemize}
Especially, the Kohn--Nirenberg transform $R$ is normalized.
\end{thm}

\begin{rem}{\rm
Condition (a) in the previous Theorem is relevant only for compact groups $G$.
}\end{rem}

\paragraph{Proof.}
Conditions (a), (b) and (d) are equivalent, because
\begin{eqnarray*}
  \iint D(u,v)(x,\eta)\,{\rm d}\eta\,{\rm d}x
  & = & FD(u,v)(\varepsilon,e) \\
  & = & \phi_D(\varepsilon,e)\,FR(u,v)(\varepsilon,e) \\
  & = & \phi_D(\varepsilon,e)\,\langle u,v\rangle.
\end{eqnarray*}
Let $a\in\mathscr S(G\times\widehat{G})$.
By Lemma~\ref{LEM:C2R},
we see that $a^D=b^R$, where
$b:=F^{-1}(\phi_D^\ast Fa)\in\mathscr S(G\times\widehat{G})$.
Hence
$$
  \iint b(x,\eta)\,{\rm d}\eta\,{\rm d}x = Fb(\varepsilon,e)
  = \phi_D(\varepsilon,e)^\ast Fa(\varepsilon,e)
  =  \phi_D(\varepsilon,e)^\ast \iint a(x,\eta)\,{\rm d}\eta\,{\rm d}x.
$$
Moreover, $b(x,\eta)=\eta(x)^\ast\,(b^R\eta)(x)$,
so that
$$
  {\rm tr}(a^D) = {\rm tr}(b^R)
  = \sum_{\eta\in\widehat{G}} d_\eta \sum_{j,k=1}^{d_\eta}
        \langle b^R\eta_{jk},\eta_{jk}\rangle
  = \iint b(x,\eta)\,{\rm d}\eta\,{\rm d}x.
$$
Thus conditions (c) and (d) are equivalent.
\hfill {\bf QED}

\begin{rem}
{\rm
Let us find how the Schwartz kernel $K\in\mathscr S(G\times G)$
of $a^D$ is related to the symbol $a\in\mathscr S(G\times\widehat{G})$
in the previous proof:
\begin{eqnarray*}
  && \langle u,a^D v\rangle \\
  & = & \langle D(u,v),a\rangle \\
  & = & \iint D(u,v)(t,\eta)\, a(t,\eta)^\ast\,{\rm d}\eta\,{\rm d}t\\
  & = & \iiint \eta(y)^\ast \int \xi(t)\,\phi_D(\xi,y)
        \int \xi(x)^\ast\,u(x)\,v(xy^{-1})^\ast\,{\rm d}x\,{\rm d}\xi\,{\rm d}y
        \,a(t,\eta)^\ast\,{\rm d}\eta\,{\rm d}t \\
  & = & \int u(x)\left(\iint \eta(y)\,a(t,\eta) \iint
        \,v(xy^{-1})\,\xi(x)\,\phi_D(\xi,y)^\ast
        \,\xi(t)^\ast\,{\rm d}\xi\,{\rm d}y\,{\rm d}\eta\,{\rm d}t
        \right)^\ast {\rm d}x.
\end{eqnarray*}
Hence we obtain
$$
  K(x,z) = \iint \eta(z^{-1}x)\,a(t,\eta)
  \int\xi(t^{-1}x)\,\phi_D(\xi,z^{-1}x)^\ast
  \,{\rm d}\xi\,{\rm d}\eta\,{\rm d}t.
$$
Here naturally $\displaystyle {\rm tr}(a^D) = \int K(x,x)\,{\rm d}x$.
}
\end{rem}

\begin{dfn}
{\rm
We say that time-frequency transform $D$ has the {\it correct time margins} if
\begin{equation}\label{EQ:margintime}
  \int D(u,v)(x,\eta)\,{\rm d}\eta = u(x)\,v(x)^\ast
\end{equation}
for all $u,v\in\mathscr S(G)$ and $x\in G$.
We say that $D$-quantization is {\it correct in time} if
\begin{equation}
  a^Dv(x) = f(x)\,v(x)
\end{equation}
for all $v\in\mathscr S(G)$ and for all symbols $a$ of the time-like form
$a(x,\eta)=f(x) I$,
where $f\in\mathscr S(G)$.
}
\end{dfn}

\begin{thm}
The following conditions are equivalent:
\begin{itemize}
\item[{\rm (a)}] $D[\delta_e]=\delta_e\otimes I$. In other words,
$
  D[\delta_e](x,\eta)=\delta_e(x)\,I.
$
\item[{\rm (b)}]
Time-frequency transform $D$ has the correct time margins.
\item[{\rm (c)}]
The $D$-quantization is correct in time.
\item[{\rm (d)}]
The ambiguity kernel satisfies
$\phi_D(\xi,e)=I$ for all
$\xi\in\widehat{G}$.
\end{itemize}
\end{thm}

\paragraph{Proof.}
For any time-frequency transform $D$ we have
$$
  F(D[\delta_e])(\xi,y)
  = \phi_D(\xi,y) \int \xi(z)\,\delta_e(z)\,\delta_e(zy^{-1})^\ast\,{\rm d}z
  = \phi_D(\xi,e)\,\delta_e(y).
$$
On the other hand, if $D[\delta_e](x,\eta)=\delta_e(x)\,I$ then
$$
  F(D[\delta_e])(\xi,y) = \int \xi(x)^\ast \int \eta(y)\,\delta_e(x)
  \,{\rm d}\eta\,{\rm d}x
= \delta_e(y)\,I.
$$
Thus conditions (a) and (d) are equivalent.
By the Fourier inverse formula,
\begin{eqnarray*}
  \int D(u,v)(x,\eta)\,{\rm d}\eta
  & = & \iint \eta(y)^\ast \int \xi(x)
  \,\phi_D(\xi,y)\,FR(u,v)(\xi,y)\,{\rm d}\xi\,{\rm d}y\,{\rm d}\eta\\
  & = & \int \xi(x)
  \,\phi_D(\xi,e)\,FR(u,v)(\xi,e)\,{\rm d}\xi \\
  & = & \int \xi(x)\,\phi_D(\xi,e)\,\widehat{u\,v^\ast}(\xi)\,{\rm d}\xi,
\end{eqnarray*}
so that conditions (b) and (d) are equivalent.
Now assume condition (b),
and let $a(x,\eta)=f(x)$. Then
\begin{eqnarray*}
  \langle u,a^D v\rangle
  & = & \langle D(u,v),a\rangle\\  
  & = & \iint D(u,v)(x,\eta)\,a(x,\eta)^\ast\,{\rm d}\eta\,{\rm d}x \\
  & = & \iint D(u,v)(x,\eta)\,{\rm d}\eta\,f(x)^\ast\,{\rm d}x \\
  & = & \int u(x)\,v(x)^\ast\,f(x)^\ast\,{\rm d}x,
\end{eqnarray*}
so $a^Dv(x)=f(x)\,v(x)$.
That is, condition (b) implies (c).
Finally, assume condition (c).
Let $(h_\alpha)_\alpha$ be a bounded left approximate identity
with $0\leq h_\alpha\in\mathscr S(G)$.
By translation, it is enough to check the time margins at $x=e$:
\begin{eqnarray*}
  \int D(u,v)(e,\eta)\,{\rm d}\eta
  & = & \iint D(u,v)(t,\eta)\,\delta_e(t)\,{\rm d}\eta\,{\rm d}t \\
  & = & \lim_\alpha \langle D(u,v),h_\alpha\otimes I\rangle \\
  & = & \lim_\alpha \langle u,(h_\alpha\otimes I)^D v\rangle \\
  & = & \lim_\alpha \langle u,h_\alpha v\rangle \\
  & = & u(e)\,v(e)^\ast.
\end{eqnarray*}
This proves condition (b) of the correct margins in time.
\hfill {\bf QED}

\begin{dfn}
{\rm
We say that time-frequency transform $D$ has
the {\it correct frequency margins} if
\begin{equation}\label{EQ:marginfrequency}
  \int D(u,v)(x,\eta)\,{\rm d}x = \widehat{u}(\eta)\,\widehat{v}(\eta)^\ast
\end{equation}
for all $u,v\in\mathscr S(G)$ and $\eta\in\widehat{G}$.
As a special case $v=u$ of \eqref{EQ:marginfrequency},
matrix
$\widehat{u}(\eta)\,\widehat{u}(\eta)^\ast$
is the ``energy density'' of $u$ at frequency
$\eta\in\widehat{G}$.
We say that $D$-quantization is {\it correct in frequency} if
\begin{equation}
  b^Dv(x) = v\ast g(x),
  \quad {\rm i.e.}\quad
  \widehat{b^Dv}(\eta)=\widehat{g}(\eta)\,\widehat{v}(\eta),
\end{equation}
for all $v\in\mathscr S(G)$ and for all symbols $b$ of the frequency-like form
$b(x,\eta)=\widehat{g}(\eta)$,
where $g\in\mathscr S(G)$.
}
\end{dfn}

\begin{thm}
The following conditions are equivalent:
\begin{itemize}
\item[{\rm (a)}] $D[{\bf 1}]={\bf 1}\otimes\delta_\varepsilon I$.
In other words,
$D[{\bf 1}](x,\eta)=\delta_\varepsilon(\eta)\,I$.
\item[{\rm (b)}]
Time-frequency transform $D$ has the correct frequency margins.
\item[{\rm (c)}]
The $D$-quantization is correct in frequency.
\item[{\rm (d)}]
The ambiguity kernel satisfies
$\phi_D(\varepsilon,y)=1\in\mathbb C$ for all $y\in G$.
\end{itemize}
\end{thm}

\paragraph{Proof.}
For any time-frequency transform $D$ we have
$$
  F(D[{\bf 1}])(\xi,y)
  = \phi_D(\xi,y) \int \xi(z)\,{\rm d}z
  = \phi_D(\varepsilon,y)\,\delta_\varepsilon(\xi).
$$
On the other hand, $D[{\bf 1}](x,\eta)=\delta_\varepsilon(\eta)\,I$ gives here
$$
  F(D[{\bf 1}])(\xi,y)
  = \int \xi(x)^\ast \int \eta(y)\,\delta_\varepsilon(\eta)\,I
  \,{\rm d}\eta\,{\rm d}x
  = \int\xi(x)^\ast\,{\rm d}x
  = \delta_\varepsilon(\xi).
$$
Hence conditions (a) and (d) are equivalent.
By $\mathscr F,\mathscr F^{-1}$ canceling each other, we obtain
\begin{eqnarray*}
  \int D(u,v)(x,\eta)\,{\rm d}x
  & = & \iint \eta(y)^\ast \int \xi(x)
  \,\phi_D(\xi,y)\,FR(u,v)(\xi,y)\,{\rm d}\xi\,{\rm d}y\,{\rm d}x\\
  & = & \int \eta(y)^\ast
  \,\phi_D(\varepsilon,y)\,FR(u,v)(\varepsilon,y)\,{\rm d}y\\
  & = & \iint \eta(y)^\ast\,\phi_D(\varepsilon,y)\,u(z)\,v(zy^{-1})^\ast
        \,{\rm d}y\,{\rm d}z.
\end{eqnarray*}
Especially,
$$
  \int D(u,\delta_e)(x,\eta)\,{\rm d}x
  = \int \eta(y)^\ast\,\phi_D(\varepsilon,y)\,u(y)\,{\rm d}y
$$
which equals to $\widehat{u}(\eta)$ for all $u\in\mathscr S(G)$
if and only if $\phi_D(\varepsilon,y)=1$ for all $y\in G$:
in that case also
$$
  D(u,v)(x,\eta)\,{\rm d}x
  = \iint \eta(z)^\ast\,u(z)\,\eta(zy^{-1})\,v(zy^{-1})^\ast
  \,{\rm d}y\,{\rm d}z
  = \widehat{u}(\eta)\,\widehat{v}(\eta)^\ast.
$$
Thus conditions (b) and (d) are equivalent.
Now assume condition (b), and let $b(x,\eta)=\widehat{g}(\eta)$.
Then
\begin{eqnarray*}
  \langle u,b^Dv\rangle
  & = & \langle D(u,v),b \rangle \\
  & = & \iint D(u,v)(x,\eta)\,b(x,\eta)^\ast\,{\rm d}\eta\,{\rm d}x \\
  & = & \iint D(u,v)(x,\eta)\,{\rm d}x\,\widehat{g}(\eta)^\ast\,{\rm d}\eta \\
  & = & \int \widehat{u}(\eta)\,\widehat{v}(\eta)^\ast\,\widehat{g}(\eta)^\ast
  \,{\rm d}\eta\\
  & = & \int \widehat{u}(\eta)
        \left(\widehat{g}(\eta)\,\widehat{v}(\eta)\right)^\ast{\rm d}\eta\\
  & = & \langle \widehat{u},\widehat{g}\,\widehat{v}\rangle
  \quad =\quad \langle u,v\ast g\rangle.
\end{eqnarray*}
Hence condition (c) follows from (b).
Finally, assume condition (c). Then
\begin{eqnarray*}
  \int D(u,v)(x,\eta)\,{\rm d}x
  & = & \iint D(u,v)(x,\omega)\,\delta_\eta(\omega)\,{\rm d}\omega\,{\rm d}x \\
  & = & \langle D(u,v),{\bf 1}\otimes\delta_\eta I\rangle \\
  & = & \langle u,({\bf 1}\otimes\delta_\eta I)^D v \rangle \\
  & = & \langle \widehat{u},\delta_\eta\,\widehat{v} \rangle \\
  & = & \widehat{u}(\eta)\,\widehat{v}(\eta)^\ast,
\end{eqnarray*}
so that we obtain condition (b) of the correct margins in frequency.
\hfill {\bf QED}

\begin{exa}
{\rm
In a sense, on a finite group $G$ of $|G|$ elements,
the minimal time-frequency transform $D$
having the correct margins would satisfy
\begin{equation}\label{EQ:addmargins}
  \phi_D(\xi,y) = \begin{cases}
    I
    & {\rm if}\ \xi=\varepsilon\ {\rm or}\ y=e,\\
    0
    & {\rm otherwise}.
    \end{cases}
\end{equation}
Then
\begin{equation}
  D(u,v)(x,\eta)
  = \widehat{u}(\eta)\,\widehat{v}(\eta)^\ast + \frac{1}{|G|}\left(
    u(x)\,v(x)^\ast-\langle u,v\rangle\right)I.
\end{equation}
Such $D$ could be added to other time-frequency transforms
that would otherwise have zero margins:
for instance, this happens
when the original localization comes from a commutator
of position and momentum operators, like on cyclic groups
in Section~\ref{SEC:cyclic}.
}
\end{exa}

\section{Positivity}

From the application point of view,
a reasonable time-frequency transform ought to be at least
normalized: this does not pose any problems.
However, it turns out that
the pointwise positivity is typically conflicting with the margin properties,
and thus positivity may not be an utterly desirable property.

\begin{dfn}
{\rm
{\it Positivity} of time-frequency transform $D$
means
$$
  D[u](x,\eta)\geq 0
$$
for all $u\in\mathscr S(G)$
and all $(x,\eta)\in G\times\widehat{G}$.
{\it Positivity} of the $D$-quantization $a\mapsto a^D$
means that for all $u\in\mathscr S(G)$
$$
  \langle u,a^Du\rangle \geq 0
$$
whenever
$a\in\mathscr S(G\times G)$ is positive in the sense that $a(x,\eta)\geq 0$
for all $(x,\eta)\in G\times\widehat{G}$.
}
\end{dfn}

\begin{exa}
{\rm
In the trivial case of the one-element group $G=\{e\}$,
defining $D(u,v)(x,\eta):=u(e)\,v(e)^\ast$
gives a positive time-frequency transform
with the correct margins in time and in frequency.
For time-frequency transforms,
positivity is a special case of symmetry:
}
\end{exa}

\begin{thm}
The following conditions are equivalent:
\begin{itemize}
\item[{\rm (a)}]
For all $u\in\mathscr S(G)$ we have $D[u](e,\varepsilon)\geq 0$.
\item[{\rm (b)}]
Time-frequency transform $D$ is positive.
\item[{\rm (c)}]
The $D$-quantization is positive.
\item[{\rm (d)}]
The time-lag kernel satisfies
$\varphi_D(x,y)=\int\kappa(x,z)\,\kappa(yx,z)^\ast {\rm d}z$
for some $\kappa$.
\end{itemize}
\end{thm}

\paragraph{Proof.}
Condition (b) trivially implies (a).
Assume condition (a).
Let $K_{\delta^D}$ denote the Schwartz kernel
of the original localization $\delta^D:\mathscr S(G)\to\mathscr S'(G)$.
Then for any $u\in\mathscr S(G)$
and $z=(z_k)_{k=1}^{d_\eta}\in\mathbb C^{d_\eta}$ we have
\begin{eqnarray*}
  \langle D[u](t,\eta)\,z,z\rangle
  & = & \sum_{j,k=1}^{d_\eta} z_j^\ast z_k\ D[u]_{jk}(t,\eta) \\
  & = & \sum_{j,k=1}^{d_\eta}
        z_j^\ast z_k \iint K_{\delta^D}(x,y)\,u(tx)\,u(ty)^\ast\,\eta_{jk}(yx^{-1})
        \,{\rm d}x\,{\rm d}y\\
  & = & \sum_{\ell=1}^{d_\eta} \iint K_{\delta^D}(x,y)
        \,u_\ell(x)\,u_\ell(y)^\ast
        {\rm d}x\,{\rm d}y \\
  & = & \sum_{\ell=1}^{d_\eta} D[u_\ell](e,\varepsilon)\ \geq\ 0,
\end{eqnarray*}
where
$\displaystyle u_\ell(x):=\sum_{k=1}^{d_\eta} z_k\,\eta_{k\ell}(x)^\ast u(tx)$.
Hence condition (a) implies (b).
Let $a\geq 0$.
Then $a^\ast=a=(a^{1/2})^2$,
where $a^{1/2}(x,\eta)=a(x,\eta)^{1/2}$
is the positive square root of $a(x,\eta)$,
and
\begin{eqnarray*}
  \langle u,a^D u\rangle
  & = & \langle D[u],a\rangle \\
  & = & \iint D[u](x,\eta)\,a(x,\eta)
  \,{\rm d}\eta\,{\rm d}x \\
  & = & \iint D[u](x,\eta)\,a(x,\eta)^{1/2}\,a(x,\eta)^{1/2}
        \,{\rm d}\eta\,{\rm d}x \\
  & = & \iint a(x,\eta)^{1/2}\,D[u](x,\eta)\,a(x,\eta)^{1/2}
        \,{\rm d}\eta\,{\rm d}x
  \quad \geq\quad 0,
\end{eqnarray*}
where the last inequality follows
because the ``integrand'' $a^{1/2}\,D[u]\,a^{1/2}$ is positive:
notice that here both the Haar integral
and the ``non-commutative $\eta$-integral'' are positive functionals.
Hence condition (b) implies (c).
Now suppose $a\mapsto a^D$ is positive and $u\in\mathscr S(G)$.
Take $(h_\alpha)_{\alpha}$ be a bounded left approximate identity,
where $0\leq h_\alpha\in\mathscr S(G)$ such that
$\lim_\alpha\langle u,h_\alpha\rangle = u(e)$.
Define $a_\alpha\in\mathscr S(G\times\widehat{G})$ by
$a_\alpha(x,\eta):=h_\alpha(x)\,\delta_\varepsilon(\eta)I$.
Then $a_\alpha(x,\eta)\geq 0$,
and
$$
  D[u](e,\varepsilon)
  = \langle D[u],\delta_{(e,\varepsilon)}\rangle
  = \lim_{\alpha}\langle D(u,u),a_\alpha\rangle
  = \lim_{\alpha}\langle u,(a_\alpha)^D u\rangle,
$$
which is non-negative due to the positivity of the quantization.
Hence condition (c) implies (a).
Assuming (d), from \eqref{EQ:original} we obtain
\begin{eqnarray*}
  D[u](e,\varepsilon)
  & \stackrel{\eqref{EQ:original}}{=} & \int u(x) \left(
    \int \varphi_D(x^{-1},y^{-1}x)^\ast\,u(y)\,{\rm d}y\right)^\ast {\rm d}x \\
  & = & \iint u(x)\,u(y)^\ast\,\varphi_D(x^{-1},y^{-1}x)\,{\rm d}y\,{\rm d}x \\
  & \stackrel{{\rm (d)}}{=} & \iiint u(x)\,u(y)^\ast
  \kappa(x^{-1},z)\,\kappa(y^{-1},z)^\ast\,{\rm d}z\,{\rm d}y\,{\rm d}x \\
  & = & \int \left| \int u(x)\,\kappa(x^{-1},z)\,{\rm d}x\right|^2{\rm d}z
  \quad \geq\quad 0,
\end{eqnarray*}
yielding condition (a).
Finally, let $\delta^D=A^2$ for a positive operator $A$. Then
\begin{eqnarray*}
  \varphi_D(x,y)
  & = & K_{\delta^D}(x^{-1},x^{-1}y^{-1})^\ast \\
  & = & \int K_A(x^{-1},z)^\ast\,K_A(z,x^{-1}y^{-1})^\ast\,{\rm d}z \\
  & = & \int K_A(x^{-1},z)^\ast\,K_A(x^{-1}y^{-1},z)\,{\rm d}z \\
  & = & \int \kappa(x,z)\,\kappa(yx,z)^\ast\,{\rm d}z,
\end{eqnarray*}
when setting $\kappa(x,z)=K_A(x^{-1},z)^\ast$.
Thus condition (d) follows from (a).
\hfill {\bf QED}

\paragraph{Spectrograms.}
A simple example of positive original localization operators
is an orthogonal projection $\delta^D:L^2(G)\to L^2(G)$
onto the $1$-dimensional subspace spanned
by a unit-energy {\it window} $w\in\mathscr S(G)$:
\begin{equation}
  \delta^D v := \langle v,w\rangle\,w.
\end{equation}
The window here should be ``focused at
$(e,\varepsilon)\in G\times\widehat{G}$'' in a reasonable sense:
most of energy of $w$ should be nearby $e\in G$,
and most of energy of $\widehat{w}$ should be nearby
$\varepsilon\in\widehat{G}$. In any case, now
$K_{\delta^D}(x,y)=w(y)^\ast w(x)$, and
\begin{eqnarray}
  D(u,v)(x,\eta)
  & = & \iint u(xz)\,\eta(z)^\ast\,K_{\delta^D}(z,y)^\ast\,\eta(y)\,v(xy)^\ast
        \,{\rm d}y\,{\rm d}z\\
  & = & \mathscr G_wu(x,\eta)\ \mathscr G_wv(x,\eta)^\ast,
\end{eqnarray}
where
\begin{equation}
  \mathscr G_wu(x,\eta) := \int \eta(y)^\ast\,u(y)\,w(x^{-1}y)^\ast\,{\rm d}y
\end{equation}
defines the $w$-windowed {\it short-time Fourier transform}
$\mathscr G_w u$ of signal $u$.
Notice that 
\begin{eqnarray}
  \mathscr G_w\delta_e(x,\eta)
  & = & w(x^{-1})^\ast\ =:\ \widetilde{w}(x), \\
  \mathscr G_w{\bf 1}(x,\eta)
  & = & \widehat{\overline{w}}(\eta)^\ast\,\eta(x)^\ast.
\end{eqnarray}
Clearly $D[u](x,\eta):=D(u,u)(x,\eta)\geq 0$,
and we may call it the {\it $w$-spectrogram} of signal $u$ at
$(x,\eta)\in G\times\widehat{G}$.
Actually, such a short-time Fourier transform formula on unimodular groups
was briefly mentioned in \cite{ChirikjianKyatkin},
as an analogue to the Euclidean case.
Let us find the corresponding ambiguity kernel $\phi_D$:
\begin{eqnarray*}
  & & FD(u,v)(\xi,y)\\
  & = & \int \xi(x)^\ast \int \eta(y) \int \eta(t)^\ast\, u(t)\, w(x^{-1}t)^\ast
        \,{\rm d}t \int w(x^{-1}s)\, v(s)^\ast\, \eta(s)\, {\rm d}s
        \,{\rm d}\eta\,{\rm d}x \\
  & = & \int \xi(x)^\ast \int u(t)\, w(x^{-1}t)^\ast
        \,w(x^{-1}ty^{-1})\, v(ty^{-1})^\ast
        \,{\rm d}t \,{\rm d}x \\
  & = & \int \left( \int \xi(x^{-1}t)\,w(x^{-1}t)^\ast\,w(x^{-1}ty^{-1})
        \,{\rm d}x \right) \xi(t)^\ast\,u(t)\,v(ty^{-1})^\ast
        \,{\rm d}t \\
  & = & \int \left(\int \xi(z)\,w(z)^\ast\,w(zy^{-1})
        \,{\rm d}z \right) \xi(t)^\ast\,u(t)\,v(ty^{-1})^\ast
        \,{\rm d}t \\
  & = & \left(\int \xi(z)^\ast\,w(z)\,w(zy^{-1})^\ast
        \,{\rm d}z \right)^\ast \int  \xi(t)^\ast\,u(t)\,v(ty^{-1})^\ast
        \,{\rm d}t \\
  & = & \phi_D(\xi,y)\,FR(u,v)(\xi,y),
\end{eqnarray*}
where
$$
  \phi_D(\xi,y)
  = \int \xi(z)\,w(z)^\ast\,w(zy^{-1})\,{\rm d}z
  = FR(w,w)(\xi,y)^\ast.
$$
Hence $\varphi_D(x,y)=(\mathscr F^{-1}\otimes I)\phi_D(x,y)
=w(x^{-1})^\ast\,w(x^{-1}y^{-1})=\widetilde{w}(x)\,\widetilde{w}(yx)^\ast$,
where $\widetilde{w}(t)=w(t^{-1})^\ast$.
The energy normalization means then
the energy normalization of the window:
$$
  1 = \phi_D(\varepsilon,e) = \int |w(x)|^2\,{\rm d}x = \|w\|^2.
$$
The correct margins in time would mean
$$
  I = \phi_D(\xi,e) = \int \xi(x)\,|w(x)|^2\,{\rm d}x
  = \widehat{|w|^2}(\xi),
$$
i.e. $|w|^2=\delta_e$, the Dirac delta at $e\in G$.
From another point of view, here
\begin{eqnarray*}
  D[\delta_e](x,\eta) & = & |w(x^{-1})|^2\ =\ |\widetilde{w}(x)|^2,\\
  D[{\bf 1}](x,\eta)
  & = & \widehat{\overline{w}}(\eta)^\ast\ \widehat{\overline{w}}(\eta).
\end{eqnarray*}
Consequently, it is too much to ask for the correct margins here,
but the energy normalization follows just from $\|w\|=1$.

\begin{rem}
{\rm
Let $D$ be a positive time-frequency transform
satisfying the correct margins both in time \eqref{EQ:margintime}
and in frequency \eqref{EQ:marginfrequency}.
Suppose $\delta^D$ is bounded on $L^2(G)$.
By the spectral decomposition of $\delta^D$,
then $G$ must be the trivial group of just one element $e$,
and $D(u,v)(x,\eta)=u(e)\,v(e)^\ast$.
}
\end{rem}

\section{Unitarity}

\begin{dfn}
{\rm
Time-frequency transform $D$ is called {\it unitary}
if it satisfies the {\it Moyal identity}
\begin{equation}\label{EQ:Moyal}
  \langle D(u,v),D(f,g)\rangle = \langle u,f\rangle\,\langle v,g\rangle^\ast
\end{equation}
for all $u,v\in\mathscr S(G)$ and $f,g\in\mathscr S'(G)$.
The $D$-quantization $a\mapsto a^D$
is called {\it unitary} if
\begin{equation}
  \langle a,b\rangle = \langle a^D,b^D\rangle
\end{equation}
for all $a,b\in\mathscr S(G\times\widehat{G})$,
where $\langle a^D,b^D\rangle = {\rm tr}\left( a^D\,(b^D)^\ast\right)$.
}
\end{dfn}

\begin{thm}\label{THM:unitary}
The following conditions are equivalent:
\begin{enumerate}
\item[{\rm (a)}]
$\langle D(u,{\bf 1}),D(\delta_e,\delta_y)\rangle = u(e)$
for all $u\in\mathscr S(G)$ and $y\in G$.
\item[{\rm (b)}]
Time-frequency transform $D$ is unitary.
\item[{\rm (c)}] The $D$-quantization is unitary.
\item[{\rm (d)}] Ambiguity operators $\phi_D(\xi,y)$ are unitary
for all $(\xi,y)\in\widehat{G}\times G$.
\end{enumerate}
Especially, the Kohn--Nirenberg transform is unitary.
\end{thm}

\begin{rem}
{\rm
In condition (a) of Theorem~\ref{THM:unitary},
on non-compact $G$ we may approximate the constant ${\bf 1}\not\in\mathscr S(G)$
within $\mathscr S(G)$.
}
\end{rem}

\paragraph{Proof.}
As $\phi_R(\xi,y)\equiv I$,
the Kohn--Nirenberg transform satisfies condition (d).
Moreover, it is unitary, because
\begin{eqnarray*}
  \langle R(u,v),R(f,g)\rangle
  & = & \iint u(x)\,\eta(x)^\ast\,\widehat{v}(\eta)^\ast
        \,\widehat{g}(\eta)\,\eta(x)\,f(x)^\ast\,{\rm d}\eta\,{\rm d}x \\
  & = & \int u(x)\,f(x)^\ast\,{\rm d}x
        \int \widehat{v}(\eta)^\ast\,\widehat{g}(\eta)
        \,{\rm d}\eta \\
  & = & \langle u,f\rangle\,\langle\widehat{g},\widehat{v}\rangle
        \quad = \quad \langle u,f\rangle\,\langle v,g\rangle^\ast.
\end{eqnarray*}
Assume (d), i.e. the unitarity of the ambiguity operators $\phi_D(\xi,y)$. Then
\begin{eqnarray*}
  && \langle D(u,v),D(f,g)\rangle \\
  & = & \langle FD(u,v),FD(f,g)\rangle \\
  & = & \iint \phi_D(\xi,y)\,FR(u,v)(\xi,y)
        \,FR(f,g)(\xi,y)^\ast\,\phi_D(\xi,y)^\ast\,{\rm d}\xi\,{\rm d}y \\
  & = & \iint FR(u,v)(\xi,y)
        \,FR(f,g)(\xi,y)^\ast
        \,{\rm d}\xi\,{\rm d}y \\
  & = & \langle FR(u,v),FR(f,g)\rangle \\
  & = & \langle R(u,v),R(f,g)\rangle.
\end{eqnarray*}
Thus condition (d) implies (b),
as we already know that $R$ is unitary.
Condition (b) implies condition (a),
because for $(u,v,f,g)=(u,{\bf 1},\delta_e,\delta_y)$
we have
$$
  \langle u,f\rangle \langle v,g\rangle^\ast
  = u(e).
$$
Now assume condition (a),
and let $(u,v,f,g)=(u,{\bf 1},\delta_e,\delta_y)$,
and $M(\omega,t):=\phi_D(\omega,t)^\ast \phi_D(\omega,t)$. Then
\begin{eqnarray*}
  u(e) & = & \langle D(u,v),D(f,g) \rangle \\
  & = & \langle FD(u,v),FD(f,g) \rangle \\
  & = & \iint M(\xi,t) \int \xi(x)^\ast u(x)\,v(xt^{-1})^\ast {\rm d}x
        \left(\int \xi(z)^\ast f(z)\,g(zt^{-1})^\ast\,{\rm d}z\right)^\ast
        {\rm d}\xi\,{\rm d}t \\
  & = & \int M(\xi,y^{-1})\,\widehat{u}(\xi)\,{\rm d}\xi.
\end{eqnarray*}
Since this holds for every $u\in\mathscr S(G)$,
we have $M(\xi,y^{-1})=I$
for every $(\xi,y)\in\widehat{G}\times G$.
Hence condition (d) follows from (a).

Finally, let us consider the Hilbert--Schmidt inner product of operators:
\begin{eqnarray*}
  \langle a^D,b^D\rangle
  & = & \iint K_{a^D}(x,y)\,K_{b^D}(x,y)^\ast\,{\rm d}x\,{\rm d}y \\
  & = & \iiint \xi(yt)\,\phi_D(\xi,t)^\ast 
  Fa(\xi,t)\,{\rm d}\xi
        \int Fb(\omega,t)^\ast\,\phi_D(\omega,t)\,\omega(yt)^\ast
        {\rm d}\omega\,{\rm d}x\,{\rm d}y \\
  & = & \iint \phi_D(\xi,t)^\ast 
  Fa(\xi,t)\,Fb(\xi,t)^\ast\,\phi_D(\xi,t)
         \,{\rm d}\xi\,{\rm d}t.
\end{eqnarray*}
It is clear that this
equals to $\langle Fa,Fb\rangle=\langle a,b\rangle$
for all $a,b\in\mathscr S(G\times\widehat{G})$
if and only if condition (d) holds:
thus conditions (c) and (d) are equivalent.
\hfill {\bf QED}

\begin{rem}
{\rm
By the previous Theorem,
unitary time-frequency transforms satisfy
the Moyal identity~\eqref{EQ:Moyal}
also for all $u,v,f,g\in L^2(G)$.
As a consequency of the unitarity of the Kohn--Nirenberg transform,
the energy densities $D[v_\alpha]$
uniformly cover the time-frequency plane $G\times\widehat{G}$
for any time-frequency transform $D$:
}
\end{rem}

\begin{cor}
Let $D$ be normalized, i.e.
$
  \phi_D(\varepsilon,e) = 1.
$
Let $(v_\alpha)_{\alpha\in J}$ be an orthonormal basis of $L^2(G)$.
Then $b^R=I$, where
\begin{equation}
  b = \sum_{\alpha\in J} D[v_\alpha].
\end{equation}
\end{cor}

\paragraph{Proof.}
Notice that
$$
  \langle u,v\rangle
  = \langle\sum_{\alpha\in J} \langle u,v_\alpha\rangle\,v_\alpha,v\rangle
  = \sum_{\alpha\in J}\langle u,v_\alpha\rangle\,\langle v,v_\alpha\rangle^\ast.
$$
Thus by the previous Theorem, for the Kohn--Nirenberg transform $R$
we have
$$
  \langle u,v\rangle
  = \sum_{\alpha\in J} \langle R(u,v),R(v_\alpha,v_\alpha)\rangle
  = \langle R(u,v),\sum_{\alpha\in J} R[v_\alpha]\rangle
  = \langle u,a^{R}v\rangle,
$$
yielding $a^{R}=I$ with
$$
  a = \sum_{\alpha\in J} R[v_\alpha].
$$
Now
$$
  \sum_{\alpha\in J} D[v_\alpha]
  = \sum_{\alpha\in J} R[v_\alpha]\ast \psi_D
  = I\ast\psi_D = \lambda I,
$$
where $\displaystyle\lambda=\iint \psi_D(x,\eta)\,{\rm d}\eta\,{\rm d}x
= \phi_D(\varepsilon,e) = 1$.
\hfill {\bf QED}

\section{Inner invariance}

Let us study the invariance under inner automorphisms
$(x\mapsto z^{-1}xz):G\to G$.
We denote $u_z(x):=u(z^{-1}xz)$ for $u\in\mathscr S(G)$
and $x,z\in G$.

\begin{dfn}
{\rm
Time-frequency transform $D$ is called {\it inner}
if it satisfies
\begin{equation}\label{EQ:inner}
  D(u_z,v_z)(x,\eta)=\eta(z)\,D(u,v)(z^{-1}xz,\eta)\,\eta(z)^\ast
\end{equation}
for all $u,v\in\mathscr S(G)$, $(x,\eta)\in G\times\widehat{G}$
and $z\in G$.
The $D$-quantization $a\mapsto a^D$
is called {\it inner} if
\begin{equation}
  \left(a^D(v_z)\right)_{z^{-1}} = a^D v
\end{equation}
for all $v\in\mathscr S(G)$ and $z\in G$
whenever
$a\in\mathscr S(G\times\widehat{G})$ satisfies
$a(z^{-1}xz,\eta)=\eta(z)^\ast\,a(x,\eta)\,\eta(z)$
for all $(x,\eta)\in G\times\widehat{G}$ and $z\in G$.
}
\end{dfn}

\begin{thm}
The following conditions are equivalent:
\begin{enumerate}
\item[{\rm (a)}]
$D[u_z](e,\varepsilon)=D[u](e,\varepsilon)$
for all $u\in\mathscr S(G)$ and $z\in G$.
\item[{\rm (b)}]
Time-frequency transform $D$ is inner.
\item[{\rm (c)}] The $D$-quantization is inner.
\item[{\rm (d)}]
$\phi_D(\xi,zyz^{-1})=\xi(z)\,\phi_D(\xi,y)\,\xi(z)^\ast$
for all $(\xi,y)\in\widehat{G}\times G$ and $z\in G$.
\end{enumerate}
Especially, the Kohn--Nirenberg transform is inner.
\end{thm}

\paragraph{Proof.}
Condition (a) is a special case of condition (b).
Condition (d) implies condition (b), because
\begin{eqnarray*}
  && D(u_z,v_z)(x,\eta) \\
  & = & \int \eta(y)^\ast \int \xi(x)\, \phi_D(\xi,y)
        \int \xi(t)^\ast\,u_z(t)\,v_z(ty^{-1})^\ast
        \,{\rm d}t\,{\rm d}\xi\,{\rm d}y \\
  & = & \int \eta(y)^\ast \int \xi(x)\, \phi_D(\xi,y)
        \int \xi(ztz^{-1})^\ast\,u(t)\,v(tz^{-1}y^{-1}z)^\ast
        \,{\rm d}t\,{\rm d}\xi\,{\rm d}y \\
  & = & \int \eta(zyz^{-1})^\ast \int \xi(z^{-1}x)\, \phi_D(\xi,zyz^{-1})
        \,\xi(z^{-1})^\ast \int \xi(t)^\ast\,\,u(t)\,v(ty^{-1})^\ast
        \,{\rm d}t\,{\rm d}\xi\,{\rm d}y \\
  & \stackrel{\rm (d)}{=} &
        \int \eta(zy z^{-1})^\ast \int \xi(z^{-1}xz)\,\phi_D(\xi,y)
        \int \xi(t)^\ast\,u(t)\,v(ty^{-1})^\ast\,{\rm d}t\,{\rm d}\xi\,{\rm d}y \\
  & = & \eta(z)\,D(u,v)(z^{-1}xz,\eta)\,\eta(z)^\ast.
\end{eqnarray*}
Suppose $a\in\mathscr S(G\times\widehat{G})$ is inner invariant: now
assuming condition (b), we obtain condition (c), because
$$
  \langle u,(a^D(v_z))_{z^{-1}}\rangle
  = \langle u_z,a^D(v_z)\rangle \\
  = \langle D(u_z,v_z),a\rangle \\
   \stackrel{\rm (b)}{=} \langle D(u,v),a\rangle \\
  = \langle u,a^D v\rangle.
$$
Now assume condition (c).
Let $(h_\alpha)_\alpha$ be an inner invariant
approximate identity in $\mathscr S(G)$.
Let $a_\alpha(x,\eta)=h_\alpha(x)\,\delta_\varepsilon(\eta)I
:\mathscr H_\eta\to\mathscr H_\eta$.
Then for all $u\in\mathscr S(G)$ and $z\in G$ we have
$$
  D[u](e,\varepsilon)
  = \langle u,\delta_{(e,\varepsilon)}^D u \rangle
  = \lim_\alpha \langle u,a_\alpha^D u\rangle
  \stackrel{\rm (c)}{=} \lim_\alpha \langle u,(a_\alpha^D(u_z))_{z^{-1}}\rangle
  = D[u_z](e,\varepsilon).
$$
Hence condition (c) implies condition (a).
Finally, conditions (a) and (d) are equivalent, because
for the kernel $\varphi_D=(\mathscr F^{-1}\otimes I)\phi_D$ on one hand
$$
  D[u](e,\varepsilon)
  = \iint \varphi_D(x^{-1},y)\,u(x)\,u(xy^{-1})^\ast\,{\rm d}x\,{\rm d}y,
$$
and on the other hand
\begin{eqnarray*}
  D[u_z](e,\varepsilon)
  & = & \iint \varphi_D(x^{-1},y)\,u_z(x)\,u_z(xy^{-1})^\ast
        \,{\rm d}x\,{\rm d}y \\
  & = & \iint \varphi_D(x^{-1},y)\,u(z^{-1}xz)\,u(z^{-1}xy^{-1}z)^\ast
        \,{\rm d}x\,{\rm d}y \\
  & = & \iint \varphi_D(zx^{-1}z^{-1},y)\,u(x)\,u(xz^{-1}y^{-1}z)^\ast
        \,{\rm d}x\,{\rm d}y \\
  & = & \iint \varphi_D(zx^{-1}z^{-1},zyz^{-1})\,u(x)\,u(xy^{-1})^\ast
        \,{\rm d}x\,{\rm d}y.
\end{eqnarray*}
This completes the proof.
\hfill {\bf QED}

\section{On locally compact groups}\label{SEC:LCG}

Time-frequency analysis on compact groups was presented above so
that the results turn out to have natural counterparts
on those locally compact groups that allow reasonable Fourier analysis.
We shall consider two families of such groups:
the Abelian ones,
and the type I second-countable unimodular locally groups.

\subsection{Locally compact Abelian groups}\label{SUBSEC:LCA}

For locally compact Abelian groups,
time-frequency analysis has been studied e.g. in \cite{Kutyniok},
and Kohn--Nirenberg pseudo-differential operators
have been treated in \cite{GrochenigStrohmer}.
We just have to modify the definitions a bit,
and then the results would hold as such.
In the commutative case, the frequency matrices would be just one-dimensional
scalars, which drastically simplifies many of the proofs.

What to change?
Let $G$ be a locally compact Abelian group.
Now $\widehat{G}$ is the {\it character group} of $G$,
consisting of the characters $\eta:G\to U(1)$,
i.e. continuous scalar unitary homomorphisms.
By the Pontryagin--van Kampen duality theorem,
$\widehat{G}$ is a locally compact Abelian group.
The group operation is given by the multiplication of the characters,
and the topology is the natural compact-open topology.
In the non-compact case,
we choose a positive regular group-invariant measure on $G$
to be the Haar measure:
this is unique up to a scalar multiple,
and $G$ has then infinite measure.
After this, we choose the Haar measure on $\widehat{G}$
so that the Fourier transform and the Fourier inverse transform formulas match:
\begin{eqnarray}
  \widehat{u}(\eta) = \int_G u(y)\,\eta(y)^\ast\,{\rm d}y, &&
  u(x)
  = \int_{\widehat{G}} \eta(x)\,\widehat{u}(\eta)\,{\rm d}\eta
\end{eqnarray}
for those $u\in L^1(G)$ for which $\widehat{u}\in L^1(\widehat{G})$.
Then we let the test function space to be
$\mathscr S(G)$, the Schwartz--Bruhat space on $G$.
The corresponding tempered distribution space is denoted by $\mathscr S'(G)$.

Why we did not choose Eymard's {\it Fourier algebra} $A(G)$
for a space of test functions on compact groups $G$?
Here $u\in A(G)$ has the norm $\|u\|_{A(G)}:=\|\widehat{u}\|_{L^1(\widehat{G})}$,
see \cite{Eymard}.
The Fourier algebra looks initially an inviting alternative,
especially as on the compact Abelian groups
it coincides with the {\it Feichtinger algebra}.
The Feichtinger algebra has turned out to be
a natural setting for time-frequency analysis
on locally compact Abelian groups, see e.g.
\cite{Feichtinger1,Feichtinger2,Grochenig}.
However, on non-commutative compact groups
the Kohn--Nirenberg transform would not map
$A(G)\times A(G)$ to $A(G\times\widehat{G})$,
and we would have the similar difficulties
with the Kohn--Nirenberg quantization,
which is our starting point for the time-frequency analysis on groups.
The difficulties boil down to that
the co-multiplication $\Delta$ does not necessarily map
$A(G)$ to $A(G\times G)$, as
\begin{eqnarray*}
  \|u\|_{A(G)}
  & = & \sum_{\eta\in\widehat{G}} d_\eta\,{\rm tr}(|\widehat{u}(\eta)|),\\
  \|\Delta u\|_{A(G\times G)}
  & = & \sum_{\eta\in\widehat{G}} d_\eta^2\,{\rm tr}(|\widehat{u}(\eta)|),
\end{eqnarray*}
where dimensions $d_\eta$ may grow arbitrarily large.
Of course, $d_\eta\equiv 1$ when the group is commutative,
and then $\Delta:A(G)\to A(G\times G)$ is an isometry,
and the Kohn--Nirenberg transform behaves well.

All in all, on a locally compact Abelian group $G$,
a time-frequency transform is a mapping
$$
  D:\mathscr S(G)\times\mathscr S(G)\to\mathscr S(G\times\widehat{G})
$$
such that
$$
  FD(u,v)(\xi,y) = \phi_D(\xi,y)\,FR(u,v)(\xi,y),
$$
where the ambiguity kernel $\phi_D:\widehat{G}\times G\to\mathbb C$
defines a Schwartz multiplier
$h\mapsto F^{-1}(\phi_D\,Fh)$.
Then we have the translation-modulation invariance
$$
  D[M_\xi T_y u](x,\eta) = D[u](x-y,\xi^{-1}y),
$$
where $T_yu(x) := u(x-y)$ and $M_\xi u(x) := \xi(x)\,u(x)$.

In case of the compact group $G$,
the approximate identities on $G\times\widehat{G}$
could be treated merely on $G$.
This is not enough on non-compact locally compact Abelian groups $G$,
but the modification for $G\times\widehat{G}$ is easy.
Notice that in the calculations for non-compact $G$,
distribution ${\bf 1}\not\in\mathscr S(G)$
occasionally has to be approximated by test functions.

\subsection{Type I second-countable unimodular groups}

Let $G$ be a type I second-countable unimodular locally compact group.
For background information, see e.g. \cite{Dixmier,Folland,Fuhr}.
Unimodularity of $G$ means that the left-invariant Haar measure
coincides with the right-invariant Haar measure:
briefly, it is the Haar measure of $G$.
Recall that a topological space is second-countable
when its topology has a countable base.
In our convention, topological groups are always Hausdorff spaces,
and consequently second-countable locally compact groups
are metrizable with a complete metric.
Moreover, second-countable locally compact groups
are of type I if and only if they are postliminal:
this means that for each $\eta\in\widehat{G}$
the compact linear operators $M:H_\eta\to H_\eta$
belong to the closure of $\{\widehat{u}(\eta):\,u\in L^1(G)\}$.

On such a group $G$,
the Schwartz--Bruhat space $\mathscr S(G)$
will be the test function space,
with the corresponding Schwartz--Bruhat distributions $\mathscr S'(G)$.
The time-frequency analysis results on compact groups are carried to $G$
without major changes in formulations and proofs.
The unit constant function ${\bf 1}:G\to\mathbb C$
is a distribution which does not belong to $\mathscr S(G)$ on non-compact $G$,
but it can be approximated by the test functions.

\section{Example of finite cyclic groups}
\label{SEC:cyclic}

Consider time-frequency analysis on the finite cyclic group
$G=\mathbb Z/N\mathbb Z$,
where $\widehat{G}\cong G$.
First, label spaces $G,\widehat{G}$ by functions $f:G\to\mathbb R$
and $\widehat{g}:\widehat{G}\to\mathbb R$.
Define respective {\it position and momentum operators}
$A,B:L^2(G)\to L^2(G)$ by
\begin{eqnarray}
  A u := f\,u, &&
  \widehat{Bu} := \widehat{g}\,\widehat{u}.
\end{eqnarray}
The {\it uncertainty observable} of measurement pair $(A,B)$ is
\begin{equation}
  \delta^{D_{\mathbb Z/N\mathbb Z}} := -{\rm i}2\pi[A,B]
  = -{\rm i}2\pi\left(AB-BA\right).
\end{equation}
This means
\begin{equation}
  \delta^{D_{\mathbb Z/N\mathbb Z}} v(x)
  = \int K_{\mathbb Z/N\mathbb Z}(x,y)\,v(y)\,{\rm d}y,
\end{equation}
where
\begin{equation}
  K_{\mathbb Z/N\mathbb Z}(x,y) =
    {\rm i}2\pi \left(f(y)-f(x)\right) g(x-y)
\end{equation}
corresponds to the time-lag kernel
$\varphi_{D_{\mathbb Z/N\mathbb Z}}:G\times G\to\mathbb C$,
\begin{equation}
  \varphi_{D_{\mathbb Z/N\mathbb Z}}(x,y) = K_{\mathbb Z/N\mathbb Z}(-x,-x-y)^\ast =
  {\rm i}2\pi \left(f(-x)-f(-x-y)\right) g(y)^\ast.
\end{equation}
As $D(u,v)(0,0)=\langle u,\delta^D v\rangle$,
by the time-frequency shift-invariance
\begin{equation}
  |D_{\mathbb Z/N\mathbb Z}(u,v)(x,\eta)|
  \leq 2\pi \|AB-BA\|\,\|u\|\|v\|
  \leq 4\pi\,\|f\|_{L^\infty} \|\widehat{g}\|_{L^\infty} \|u\|\,\|v\|
\end{equation}
for all $(x,\eta)\in G\times\widehat{G}$.
For the ambiguity kernel
$\phi_{D_{\mathbb Z/N\mathbb Z}}:\widehat{G}\times G\to\mathbb C$,
\begin{equation}
  \phi_{D_{\mathbb Z/N\mathbb Z}}(\xi,y) = {\rm i}2\pi
  \widehat{f}(-\xi) \left(1-{\rm e}^{{\rm i}2\pi\xi y/N}\right) g(y)^\ast.
\end{equation}
A natural choice for the position labeling function
$f:G\to\mathbb R$ could be
\begin{equation}\label{EQ:f-left}
  f(x) := x/N\quad {\rm for}\quad 0\leq x<N
\end{equation}
(here $f(x):=x/N$ for $0<x\leq N$ would be another good choice,
but it ultimately leads to the same limit as $N\to\infty$ in the next section).
Observe that for $0<\eta<N$
$$
  0
  = N^{-1}\sum_{x=0}^{N-1}\left((x+1)/N-x/N\right)
  {\rm e}^{-{\rm i}2\pi x\eta/N}
  = {\rm e}^{{\rm i}2\pi\eta/N}\left(\widehat{f}(\eta) + N^{-1}\right)
  -\widehat{f}(\eta),
$$
yielding
\begin{equation}
  \widehat{f}(\eta) = \frac{-1/N}{1-{\rm e}^{-{\rm i}2\pi\eta/N}},
\end{equation}
so that if $\widehat{g}(\eta)=f(\eta)$
(i.e. $g(y)=N\widehat{f}(-y)$) then
\begin{equation}\label{EQ:ambiguitydiscrete}
  \phi_{D_{\mathbb Z/N\mathbb Z}}(\xi,y) = \begin{cases}
    \frac{{\rm i}2\pi}{N}
    \frac{1-{\rm e}^{{\rm i}2\pi\xi y/N}}
    {\left(1-{\rm e}^{{\rm i}2\pi\xi/N}\right)
      \left(1-{\rm e}^{-{\rm i}2\pi y/N}\right)}
    & {\rm if}\ \xi\not=0\ {\rm and}\ y\not=0,\\
  0 & {\rm if}\ \xi=0\ {\rm or}\ y=0.
  \end{cases}
\end{equation}
Let us define the time-frequency transform $Q_{\mathbb Z/N\mathbb Z}$
on the finite cyclic group $G=\mathbb Z/N\mathbb Z$
by its ambiguity kernel, where
\begin{equation}
  \phi_{Q_{\mathbb Z/N\mathbb Z}}(\xi,y) = \begin{cases}
    \frac{{\rm i}2\pi}{N}
    \frac{1-{\rm e}^{{\rm i}2\pi\xi y/N}}
    {\left(1-{\rm e}^{{\rm i}2\pi\xi/N}\right)
      \left(1-{\rm e}^{-{\rm i}2\pi y/N}\right)}
    & {\rm if}\ \xi\not=0\ {\rm and}\ y\not=0,\\
  1 & {\rm if}\ \xi=0\ {\rm or}\ y=0.
  \end{cases}
\end{equation}
That is, we summed \eqref{EQ:ambiguitydiscrete} and \eqref{EQ:addmargins},
obtaining the correct margins.

\begin{thm}
Mapping $[u]\mapsto Q_{\mathbb Z/N\mathbb Z}[u]$
is invertible for all $N\in\mathbb Z^+$.
The corresponding $Q_{\mathbb Z/N\mathbb Z}$-quantization
is invertible if and only if $N$ is prime or $N=1$.
\end{thm}

\paragraph{Proof.}
Let $D=Q_{\mathbb Z/N\mathbb Z}$. Case $N=1$ is trivial.
Assume now that $N$ is prime.
Then ambiguity kernel $\phi_D$ has no zeros,
so let $g(\xi,y):=\phi_D(\xi,y)^{-1}$.
Hence starting from
$D[u] = F^{-1}(\phi_D\,FR[u])$
we find $FR[u]=g\,FD[u]$,
and from it we obtain $u(x)\,u(x-y)^\ast$ for all $x,y\in\mathbb Z/N\mathbb Z$.
Thus $[u]\mapsto D[u]$ is invertible when $N$ is prime.
What about the invertibility of the $D$-quantization $a\mapsto a^D$?
Recall that the Kohn--Nirenberg quantization $a\mapsto a^R$
is invertible:
linear mapping $A:L^2(\mathbb Z/N\mathbb Z)\to L^2(\mathbb Z/N\mathbb Z)$
is of the form $A=a^R$,
where
$
  a(x,\eta)=\eta(x)^\ast A\eta(x)
$
for $\eta(x):={\rm e}^{{\rm i}2\pi x\eta/N}$.
Then the $D$-quantization $b\mapsto b^D$
is invertible, because
\begin{eqnarray*}
  && \langle u,a^R v\rangle
  = \langle R(u,v),a\rangle
  = \langle FR(u,v),Fa\rangle
  = \langle FD(u,v),Fb\rangle
  = \langle D(u,v),b\rangle \\
  & = & \langle u,b^D v\rangle,
\end{eqnarray*}
where $Fb=g^\ast Fa$:
here $a^R=b^D$.
This concludes the case of prime $N$.

Finally, let us consider divisible $N\geq 4$.
Now $\phi_D(\xi,y)=0$ if and only if
$\xi,y$ are zero divisors modulo $N$.
In this case, $b^D=0$ if $b$ is a symbol such that $Fb$ is supported
only on the zero divisors.
Hence the $D$-quantization is not injective, nor surjective
(due to the finite-dimensionality).
However, it turns out that $[u]\mapsto D[u]$ is still invertible.
Finding $[u]$ from $D[u]$
is reduced to phase retrieval,
as we easily get the time margins
$
  |u(x)|^2 = \sum_{\eta=1}^N D[u](x,\eta).
$
Especially, case $u=0$ is trivial, so assume $u\not=0$.
Knowing $D[u]$, we also find
$$
  F^{-1}D[u](\xi,y) =
  \phi_D(\xi,y)\,\frac{1}{N}\sum_{z=1}^N {\rm e}^{-{\rm i}2\pi z\xi/N}
  \,u(z)\,u(z-y)^\ast.
$$
From this, since
$\displaystyle
\frac{1-{\rm e}^{{\rm i}2\pi\xi y/N}}{1-{\rm e}^{{\rm i}2\pi\xi/N}}
= \sum_{k=0}^{y-1} {\rm e}^{{\rm i}2\pi k\xi/N}$
for $0<y<N$,
we obtain numbers
$$
  E(x,y) := \sum_{k=0}^{y-1} u(x+k)\,u(x+k-y)^\ast
$$
for all $x$.
We may recover only the equivalence class $[u]$ of $u$,
but suppose we know the complex phase of some $u(z)\not=0$.
We proceed recursively as follows:
We find numbers $u(z+1)$ and $u(z-1)$
from $E(z+1,1)$ and $E(z,1)$, respectively.
If we have already recovered numbers
$u(z\pm h)$ for $0\leq h<j$,
then we stably obtain numbers of $u(z+j)$ and $u(z-j)$
by finding their complex phases
from $E(z+1,j)$ and $E(z,j)$, respectively.
This completes the proof.
\hfill{\bf QED}

\begin{rem}
{\rm
In the previous proof,
the stable algorithm for $D[u]\mapsto [u]$
can be built around any point $z\in\mathbb Z/N\mathbb Z$ for which $u(z)\not=0$.
Let us also note the estimates
$$
  |\phi_D(\xi,y)|\leq |\phi_D(1,y)|
  = \frac{2\pi}{N}\left|1-{\rm e}^{{\rm i}2\pi/N}\right|^{-1}\leq\frac{\pi}{2}
$$
for all $N\geq 2$ and $\xi,y$.
Without losing generality, for $0<y\leq N/2$ this follows by observing that 
$$
  \phi_D(\xi,y) = \frac{{\rm i}2\pi}{N}
  \left(1-{\rm e}^{-{\rm i}2\pi y/N}\right)^{-1}
  \sum_{k=0}^{y-1} {\rm e}^{{\rm i}2\pi k\xi/N}.
$$
By the geometry of the unit circle, the optimal bounds
$$
  |\phi_D(\xi,y)|\leq \frac{2\pi}{N} \left|1-{\rm e}^{{\rm i}2\pi/N}\right|^{-1}
$$
form a monotonically decreasing sequence with the limit $1$ as $N\to\infty$.
}
\end{rem}

\section{Limit of cyclic case: Born--Jordan}

Next we study what happens to transforms $D_{\mathbb Z/N\mathbb Z}$
when we take the limit $N\to\infty$ interpreting
either that $\mathbb Z/N\mathbb Z$ tends to the compact circle group
$\mathbb T=\mathbb R/\mathbb Z$ or to the non-compact group
$\mathbb Z$ of integers.
We also study the further limiting time-frequency transforms
on the real line $\mathbb R$.

Starting from natural time-frequency transforms
of signals on $\mathbb Z/N\mathbb Z$,
we study the limiting cases on
compact $\mathbb T$ and non-compact $\mathbb Z$,
and their limits on $\mathbb R$.
At limit $N\to\infty$ to compact group $\mathbb T$,
from transforms $D_{\mathbb Z/N\mathbb Z}$ in the previous section
we obtain time-frequency transform $D_{\mathbb T}$ with ambiguity kernel
$\phi_{D_{\mathbb T}}:\mathbb Z\times\mathbb T\to\mathbb C$,
where
\begin{equation}
  \phi_{D_{\mathbb T}}(\xi,y) = \begin{cases} -\xi^{-1}
  \left(1-{\rm e}^{{\rm i}2\pi\xi y}\right) /
    \left(1-{\rm e}^{-{\rm i}2\pi y}\right) &
    {\rm if}\ \xi\not=0\ {\rm and}\ y\not=0, \\
    1 & {\rm if}\ \xi\not=0\ {\rm and}\ y=0, \\
    0 & {\rm if}\ \xi=0.
    \end{cases}
\end{equation}
Indeed, $y\mapsto\phi_{D_{\mathbb T}}(\xi,y)$
is a trigonometric polynomial:
\begin{eqnarray*}
  \phi_{D_\mathbb T}(\xi,y)
  & = &
  \frac{1}{|\xi|}\sum_{k=0}^{|\xi|-1}{\rm e}^{-{\rm i}2\pi yk}
        \quad{\rm if}\quad \xi<0,\\
  \phi_{D_\mathbb T}(\xi,y)
  & = & \frac{1}{\xi}
  \sum_{k=1}^{\xi} {\rm e}^{+{\rm i}2\pi yk}
        \quad\quad{\rm if}\quad \xi>0.
\end{eqnarray*}
Indeed,
$(h\mapsto F^{-1}(\phi_D\,Fh)):
\mathscr S(\mathbb T\times\mathbb T)\to\mathscr S(\mathbb T\times\mathbb Z)$
is a Schwartz multiplier.
Moreover, time-frequency transform $D_{\mathbb T}$ is band-limited, mapping
$\mathscr T(\mathbb T)\times\mathscr T(\mathbb T)$
to $\mathscr T(\mathbb T\times\mathbb Z)$.
Since $|\phi_{D_{\mathbb T}}(\xi,y)|\leq 1$,
by Theorem~\ref{THM:L2boundedness} we have the $L^2$-bounds
\begin{eqnarray}
  \|D_{\mathbb T}(u,v)\|\leq \|u\|\,\|v\|,\\
  \|a^{D_\mathbb T} v\|\leq \|a\|\,\|v\|.
\end{eqnarray}
Analogously, we have time-frequency transform
$D_{\mathbb Z}$ on non-compact group $\mathbb Z$, with
ambiguity kernel
$\phi_{D_{\mathbb Z}}:\mathbb T\times\mathbb Z\to\mathbb C$,
\begin{equation}
  \phi_{D_{\mathbb Z}}(\xi,y) = \begin{cases}
    y^{-1}
  \left(1-{\rm e}^{{\rm i}2\pi\xi y}\right) /
    \left(1-{\rm e}^{{\rm i}2\pi\xi}\right) &
    {\rm if}\ \xi\not=0\ {\rm and}\ y\not=0, \\
    1 & {\rm if}\ \xi=0\ {\rm and}\ y\not=0,  \\
    0 & {\rm if}\ y=0.
    \end{cases}
\end{equation}
Hence time-lag kernel
$\varphi_{D_{\mathbb Z}}:\mathbb Z\times\mathbb Z\to\mathbb C$ is given by
\begin{equation}
  \varphi_{D_{\mathbb Z}}(x,y) = \begin{cases}
    1/|y| & {\rm if}\ -y< x\leq 0\ {\rm or}\ 0< x\leq -y,\\
    0 & {\rm otherwise.}
    \end{cases}
\end{equation}
Here $\varphi_{D_{\mathbb Z}}(x,y)=K_{\mathbb Z}(-x,-x-y)^\ast$
(equivalently, $K_{\mathbb Z}(x,y)=\varphi_{D_{\mathbb Z}}(-x,x-y)^\ast$),
with
\begin{equation}
  \delta^{D_{\mathbb Z}} v(x) = \sum_{y\in\mathbb Z} K_{\mathbb Z}(x,y)\,v(y),
\end{equation}
with kernel $K_{\mathbb Z}:\mathbb Z\times\mathbb Z\to\mathbb C$ given by
\begin{equation}
  K_{\mathbb Z}(x,y) = \begin{cases}
    1/|x-y| & {\rm if}\ y<0\leq x\ {\rm or}\ x<0\leq y,\\
    0 & {\rm otherwise.}
    \end{cases}
\end{equation}
At the continuum limit on $\mathbb R$, we obtain
time-frequency transform $D_{\mathbb R}$, with
\begin{equation}
  \delta^{D_{\mathbb R}}v(x) = \int_{\mathbb R} K_{\mathbb R}(x,y)\,v(y)\,{\rm d}y,
\end{equation}
where Schwartz kernel $K_{\mathbb R}:\mathbb R\times\mathbb R\to\mathbb C$
is given by
\begin{equation}
  K_{\mathbb R}(x,y) = \begin{cases}
    1/|x-y| & {\rm if}\ xy<0, \\
    0 & {\rm otherwise.}
    \end{cases}
\end{equation}
Hence $D_{\mathbb R}=Q$ is the Born--Jordan transform,
\begin{equation}
  Q(u,v)(x,\eta) = \int_{\mathbb R}
  {\rm e}^{-{\rm i}2\pi y\eta}\frac{1}{y}\int_{x-y/2}^{x+y/2}
  u(t+y/2)\,v(t-y/2)^\ast
  \,{\rm d}t\,{\rm d}y.
\end{equation}
Time-frequency transform $D_{\mathbb Z/N\mathbb Z}$
has zero margins in both time and in frequency,
but the margins for $Q$ are correct.

\paragraph{Alternative way.}
Above, we went from $\mathbb Z/N\mathbb Z$ to $\mathbb R$ via $\mathbb Z$.
What if our route would have been via $\mathbb T$ instead?
The outcome must still be the Born--Jordan transform.
Let us check this process:
Time-frequency transform
$D_{\mathbb T}$ on compact group $\mathbb T$ has
time-lag kernel
$\varphi_{D_{\mathbb T}}:\mathbb T\times\mathbb T\to\mathbb C$, where
for $y\not=0$ we have
\begin{equation}
  \varphi_{D_{\mathbb T}}(x,y) =
    {\rm i}2\pi\,\frac{w(x)-w(x+y)}{1-{\rm e}^{-{\rm i}2\pi y}},
\end{equation}
with the sawtooth wave $w:\mathbb T\to\mathbb R$
satisfying $w(x)=x$ for $0<x<1$.
Now
\begin{equation}
  \delta^{D_{\mathbb T}} v(x)
  = \int K_{\mathbb T}(x,y)\,v(y)\,{\rm d}y,
\end{equation}
with kernel $K_{\mathbb T}:\mathbb T\times\mathbb T\to\mathbb C$ given by
$K_{\mathbb T}(x,y)=\varphi_{D_{\mathbb T}}(-x,x-y)^\ast$,
\begin{equation}
  K_{\mathbb T}(x,y)
  = -{\rm i}2\pi\,\frac{1-(x-y)}{1-{\rm e}^{{\rm i}2\pi(x-y)}}
\end{equation}
when $-1<y<0<x<1$ and $x-y\not=1$: if here $x,y\to 0$,
we again obtain the Born--Jordan transform $Q$ as the continuum limit.
Properties of the Born--Jordan transform were studied in \cite{Turunen},
where also closely related variants of $D_{\mathbb T},D_{\mathbb Z}$
were introduced.

\section{Computed pictures of discrete distributions}

In the following pictures,
we present three different discrete time-frequency distributions
for the same signal:
the periodic and non-periodic Born--Jordan distributions, and a spectrogram.
The original speech signal of the author has 1000 samples,
with sampling rate of 4000 Hz. The pictures were produced using Matlab.
In the grey-scale time-frequency distribution pictures,
higher values are darker in shade.
For the spectrogram, zero value corresponds to white.
For the other time-frequency images, zero value corresponds to mid-grey.

\newpage
\includegraphics[trim=10mm 0mm 0mm 10mm,scale=0.75]{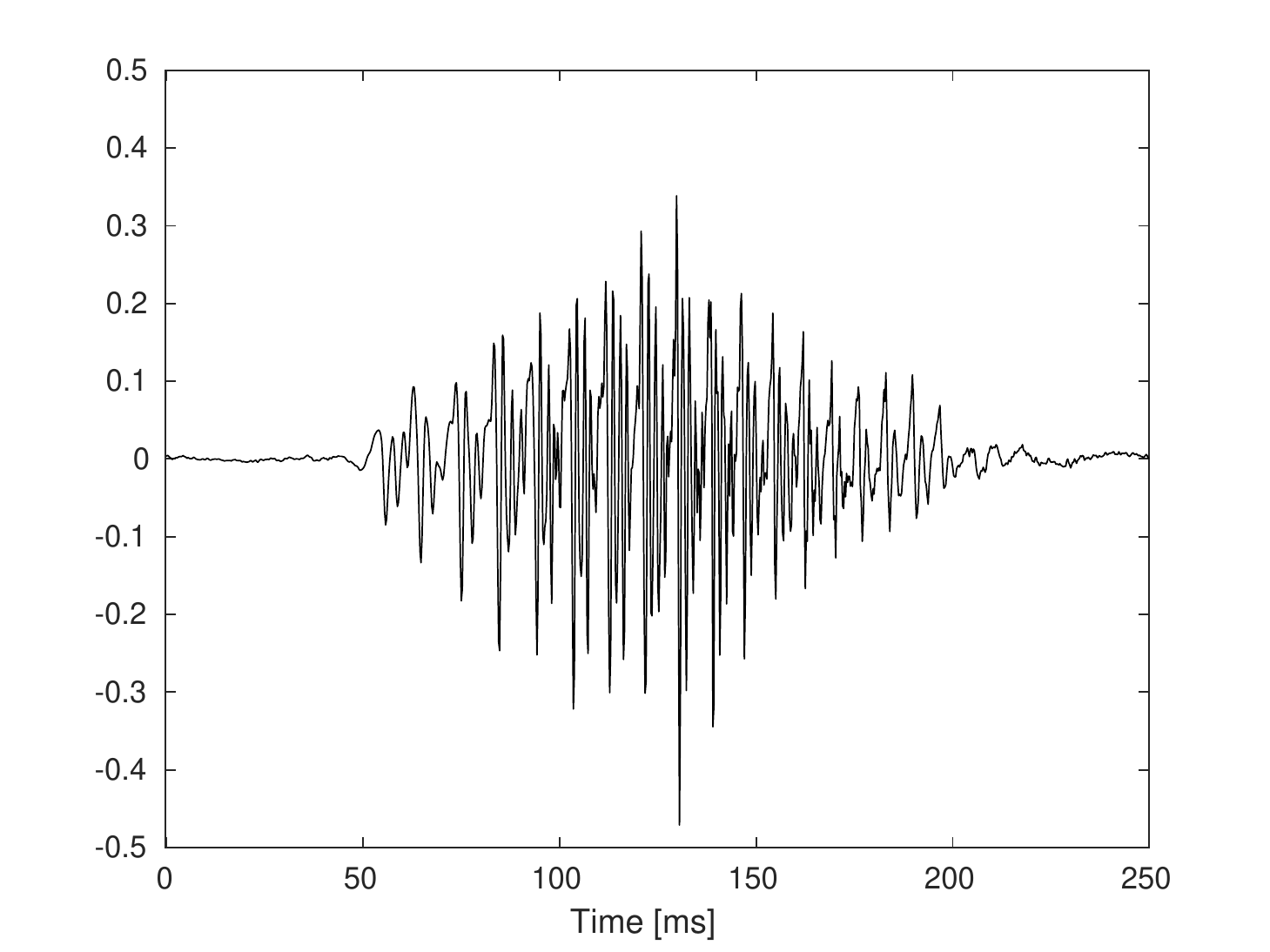}
\captionof{figure}
{Speech signal ``Why?'', sampling rate 4000 Hz.}

\includegraphics[trim=3mm 0mm 0mm 0mm,scale=0.75]{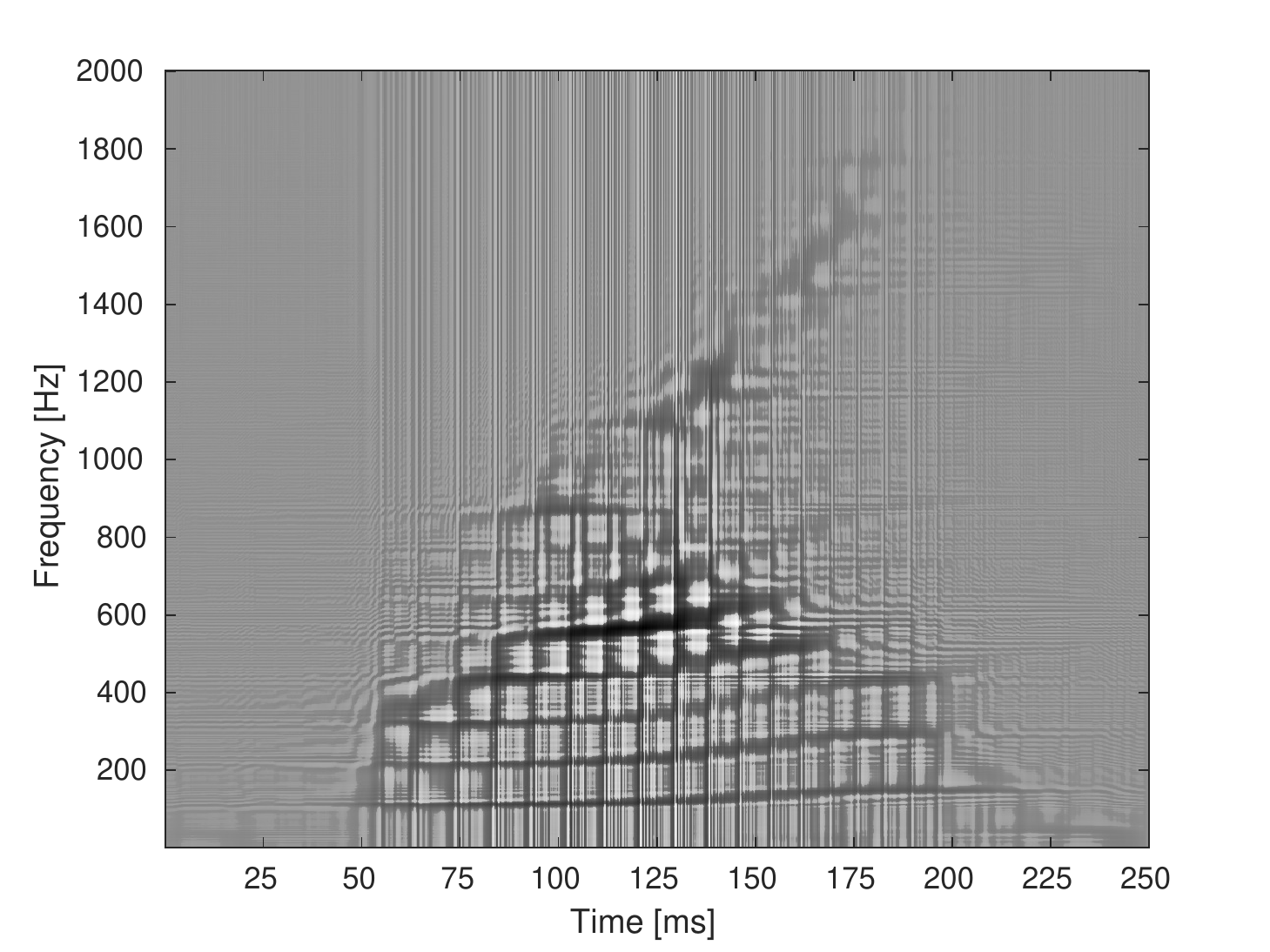}
\captionof{figure}{Time-frequency distribution $Q_{\mathbb Z}[u]$ for signal $u$ (``Why?'').}

\newpage

\includegraphics[trim=3mm 0mm 0mm 0mm,scale=0.75]{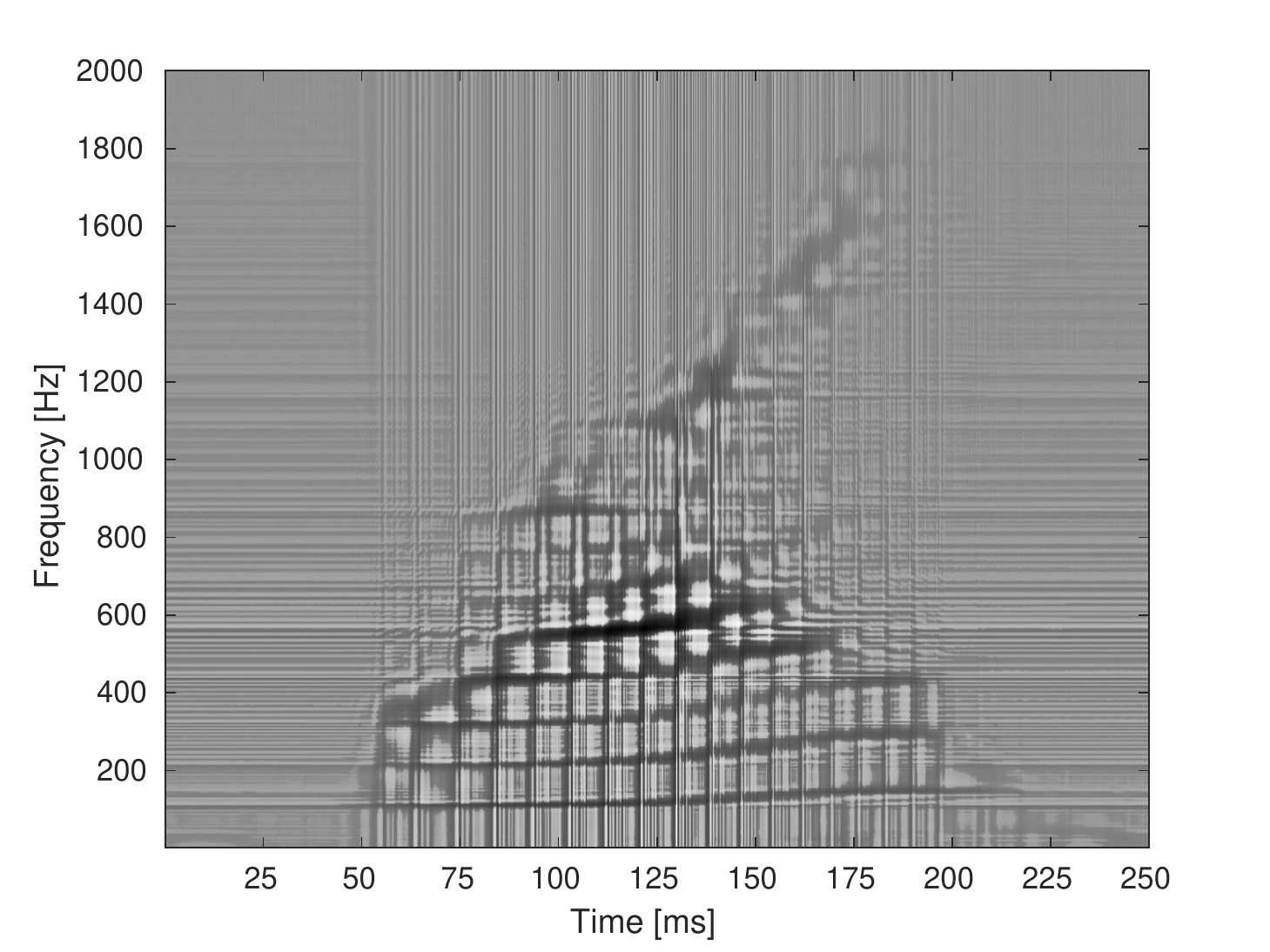}
\captionof{figure}{Time-frequency distribution
$Q_{\mathbb Z/N\mathbb Z}[u]$ for the periodized signal $u$
(``...Why\,Why\,Why\,Why...''), zooming into a single period of 250 ms.}

\includegraphics[trim=3mm 0mm 0mm 0mm,scale=0.75]{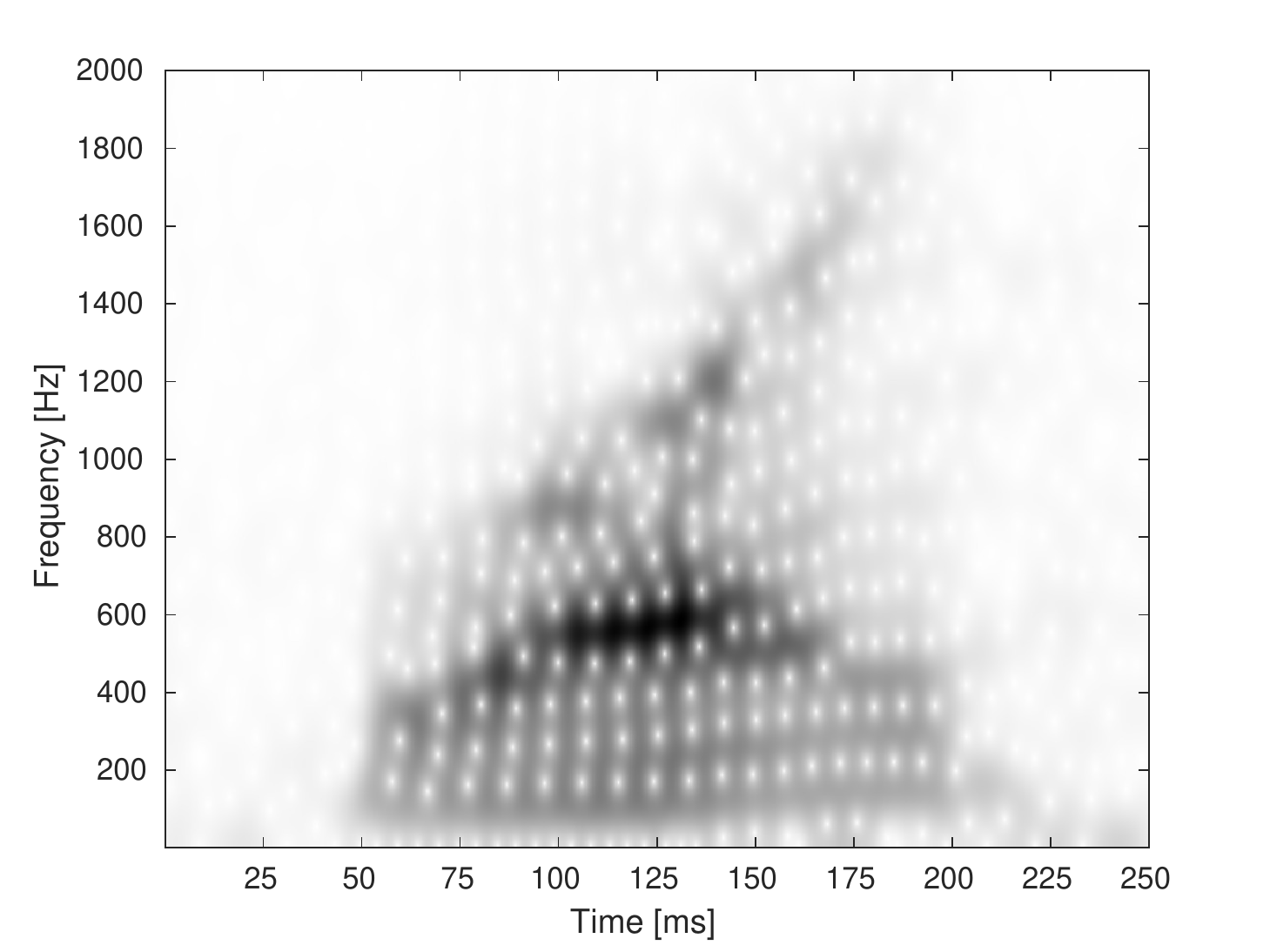}
\captionof{figure}{Spectrogram for periodized signal $u$
(``...Why\,Why\,Why\,Why...''), with a Gaussian window,
zooming into a single period of 250 ms.}

\end{document}